\documentclass[10pt]{elsart}
\usepackage{amssymb,amsmath}
\usepackage[letterpaper,height=9in]{geometry}
\usepackage{subfigure,url,warmread,graphicx}
\usepackage[all,import]{xy}

\newcommand{\norm}[1]{\ensuremath{\left\| #1 \right\|}}
\newcommand{\bracket}[1]{\ensuremath{\left[ #1 \right]}}
\newcommand{\braces}[1]{\ensuremath{\left\{ #1 \right\}}}
\newcommand{\parenth}[1]{\ensuremath{\left( #1 \right)}}
\newcommand{\refeqn}[1]{(\ref{eqn:#1})}
\newcommand{\reffig}[1]{Fig. \ref{fig:#1}}
\newcommand{\tr}[1]{\mbox{tr}\ensuremath{\negthickspace\bracket{#1}}}
\newcommand{\deriv}[2]{\ensuremath{\frac{\partial #1}{\partial #2}}}
\newcommand{\SO}{\ensuremath{\mathrm{SO(3)}}}
\newcommand{\so}{\ensuremath{\mathfrak{so}(3)}}
\newcommand{\SE}{\ensuremath{\mathrm{SE(3)}}}

\newcommand{\aSE}[2]{\ensuremath{\begin{bmatrix}#1&#2\\0&1\end{bmatrix}}}

\begin{document}
\allowdisplaybreaks

\begin{frontmatter}

\title{\LARGE Lie Group Variational Integrators for the Full Body Problem}

\author[AE]{Taeyoung Lee\thanksref{Rackham}},
\ead{tylee@umich.edu}
\author[MA]{Melvin Leok\thanksref{Rackham}}, and
\ead{mleok@umich.edu}
\author[AE]{N. Harris McClamroch\thanksref{NSF}}
\ead{nhm@engin.umich.edu}
\address[AE]{Department of Aerospace Engineering, University of Michigan, Ann Arbor, MI 48109, USA}
\address[MA]{Department of Mathematics, University of Michigan, Ann Arbor, MI 48109, USA}

\thanks[Rackham]{\renewcommand{\baselinestretch}{1.0}\footnotesize This research has been supported in part by a grant from the Rackham Graduate School, University of Michigan.}
\thanks[NSF]{This research has been supported in part by NSF under
grant ECS-0244977.}

\begin{abstract}
We develop the equations of motion for full body models that
describe the dynamics of rigid bodies, acting under their mutual
gravity. The equations are derived using a variational approach
where variations are defined on the Lie group of rigid body
configurations. Both continuous equations of motion and variational
integrators are developed in Lagrangian and Hamiltonian forms, and
the reduction from the inertial frame to a relative coordinate
system is also carried out. The Lie group variational integrators
are shown to be symplectic, to preserve conserved quantities, and to
guarantee exact evolution on the configuration space. One of these
variational integrators is used to simulate the dynamics of two
rigid dumbbell bodies.
\end{abstract}


\end{frontmatter}

{\renewcommand{\baselinestretch}{1.09}\footnotesize

\section*{Nomenclature}
\begin{tabular}{cp{12cm}c}
$\gamma_i$ & Linear momentum of the $i$th body in the inertial frame & p.\pageref{no:gammai}\\
$\Gamma$ & Relative linear momentum & p.\pageref{no:Gamma}\\
$J_i$ & Standard moment of inertia matrix of the $i$th body & p.\pageref{no:Ji}\\
$J_{d_i}$ & Nonstandard moment of inertia matrix of the $i$th body & p.\pageref{no:Jdi}\\
$J_R$ & Standard moment of inertia matrix of the first body with respect to the second body fixed frame & p.\pageref{no:JR}\\
$J_{d_R}$ & Nonstandard moment of inertia matrix of the first body with respect to the second body fixed frame & p.\pageref{no:JdR}\\
$m_i$ & Mass of the $i$th body & p.\pageref{no:mi}\\
$m$ & Reduced mass for two bodies of mass $m_1$ and $m_2$, $m=\frac{m_1m_2}{m_1+m_2}$ & p.\pageref{no:m}\\
$M_i$ & Gravity gradient moment on the $i$th body & p.\pageref{no:Mi}\\
$\Omega_i$ & Angular velocity of the $i$th body in its body fixed frame & p.\pageref{no:Omegai}\\
$\Omega$ & Angular velocity of the first body expressed in the second body fixed frame & p.\pageref{no:Omega}\\
$\mathbf{\Omega}$ & $\mathbf{\Omega}=(\Omega_1,\Omega_2,\cdots,\Omega_n)$ & p.\pageref{no:bfOmega}\\
$\Pi_i$ & Angular momentum of the $i$th body in its body fixed frame & p.\pageref{no:Pii}\\
$\Pi$ & Angular momentum of the first body expressed in the second body fixed frame & p.\pageref{no:Pi}\\
$R_i$ & Rotation matrix from the $i$th body fixed frame to the inertial frame & p.\pageref{no:Ri}\\
$R$ & Relative attitude of the first body with respect to the second body & p.\pageref{eqn:R}\\
$\mathbf{R}$ & $\mathbf{R}=(R_1,R_2,\cdots,R_n)$ & p.\pageref{no:bfR}\\
$v_i$ & Velocity of the mass center of the $i$th body in the inertial frame & p.\pageref{no:vi}\\
$V$ & Relative velocity of the first body with respect to the second body in the second body fixed frame & p.\pageref{no:V}\\
$V_2$ & Velocity of the second body in the second body fixed frame & p.\pageref{eqn:V2}\\
$x_i$ & Position of the mass center of the $i$th body in the inertial frame & p.\pageref{no:xi}\\
$X$ & Relative position of the first body with respect to the second body expressed in the second body fixed frame & p.\pageref{eqn:X}\\
$X_2$ & Position of the second body in the second body fixed frame & p.\pageref{eqn:X2}\\
$\mathbf{x}$ & $\mathbf{x}=(x_1,x_2,\cdots,x_n)$ & p.\pageref{no:bfx}\\
\end{tabular}}
\newpage

\footnotesize

\section{Introduction}
\subsection{Overview}
The full body problem studies the dynamics of rigid
bodies interacting under their mutual potential, and the mutual
potential of distributed rigid bodies depends on both the position
and the attitude of the bodies. Therefore, the translational and the
rotational dynamics are coupled in the full body problem. The full
body problem arises in numerous engineering and scientific fields.
For example, in astrodynamics, the trajectory of a large spacecraft
around the Earth is affected by the attitude of the spacecraft, and
the dynamics of a binary asteroid pair is characterized by the
non-spherical mass distributions of the bodies. In chemistry, the
full rigid body model is used to study molecular dynamics. The
importance of the full body problem is summarized
in~\cite{pro:Koon03}, along with a preliminary discussion from the
point of view of geometric mechanics.

The full two body problem was studied by
Maciejewski~\cite{jo:macie}, and he presented equations of motion in
inertial and relative coordinates and discussed the existence of
relative equilibria in the system. However, he does not derive the
equations of motion, nor does he discuss the reconstruction
equations that allow the recovery of the inertial dynamics from the
relative dynamics. Scheeres derived a stability condition for the
full two body problem~\cite{jo:Scheeres02}, and he studied the
planar stability of an ellipsoid-sphere model~\cite{pro:Scheeres03}.
Recently, interest in the full body problem has increased, as it is
estimated that up to 16\% of near-earth asteroids are
binaries~\cite{jo:margot}. Spacecraft motion about binary asteroids
have been discussed using the restricted three body
model~\cite{pro:Scheeres03b},~\cite{pro:gabern}, and the four body
model~\cite{jo:Scheeres05}.

Conservation laws are important for studying the
dynamics of the full body problem, because they describe qualitative
characteristics of the system dynamics. The representation used for
the attitude of the bodies should be globally defined since the
complicated dynamics of such systems would require frequent
coordinate changes when using representations that are only defined
locally. General numerical integration methods, such as the
Runge-Kutta scheme, do not preserve first integrals nor the exact
geometry of the full body dynamics~\cite{bk:hairer}. A more careful
analysis of computational methods used to study the full body
problem is crucial.

Variational integrators and Lie group methods provide a systematic
method of constructing structure-preserving numerical
integrators~\cite{bk:hairer}. The idea of the variational approach
is to discretize Hamilton's principle rather than the continuous
equations of motion~\cite{jo:marsden}. The numerical integrator
obtained from the discrete Hamilton's principle exhibits excellent
energy properties~\cite{jo:hairer}, conserves first integrals
associated with symmetries by a discrete version of Noether's
theorem, and preserves the symplectic structure. Many interesting
differential equations, including full body dynamics, evolve on a
Lie group. Lie group methods consist of numerical integrators that
preserve the geometry of the configuration space by automatically
remaining on the Lie group~\cite{ic:iserles}.

Moser and Vesselov~\cite{jo:moser}, Wendlandt and
Marsden~\cite{jo:wendlandt} developed numerical integrators for a
free rigid body by imposing an orthogonal constraint on the attitude
variables and by using unit quaternions, respectively. The idea of
using the Lie group structure and the exponential map to numerically
compute rigid body dynamics arises in the work of Simo, Tarnow, and
Wong~\cite{jo:simo}, and in the work by Krysl~\cite{Kr2004}. A Lie
group approach is explicitly adopted by Lee, Leok, and McClamroch in
the context of a variational integrator for rigid body attitude
dynamics with a potential dependent on the attitude, namely the 3D
pendulum dynamics, in~\cite{pro:Lee05}.

The motion of full rigid bodies depends essentially on the mutual
gravitational potential, which in turn depends only on the relative
positions and the relative attitudes of the bodies. Marsden et al.
introduce discrete Euler-Poincar\'{e} and Lie-Poisson equations
in~\cite{jo:marsden99} and \cite{jo:marsden00}. They reduce the
discrete dynamics on a Lie group to the dynamics on the
corresponding Lie algebra. Sanyal, Shen and McClamroch develop
variational integrators for mechanical systems with configuration
dependent inertia and they perform discrete Routh reduction
in~\cite{jo:sanyal}. A more intrinsic development of discrete Routh
reduction is given in~\cite{leok} and~\cite{jo:Leok}.

\subsection{Contributions}
The purpose of this paper is to provide a complete set of equations
of motion for the full body problem. The equations of motion are
categorized by three independent properties: continuous / discrete
time, inertial / relative coordinates, and Lagrangian / Hamiltonian
forms. Therefore, a total of eight types of equations of motion for
the full body problem are given in this paper. The relationships
between these equations of motion are shown in \reffig{8eom}, and
are further summarized in \reffig{cubic}.

\renewcommand{\thefigure}{\arabic{figure}}
\begin{figure}[htbp]
\resizebox{1.00\linewidth}{!}{
\setlength{\unitlength}{1.65em}\centering\footnotesize{
\begin{picture}(22.3,11)(0,-11)
\put(6.15,-1.0){\framebox(10,1)[c]
{\shortstack[c]{Continuous time (Sec. 3)}}}
\put(11.15,-1){\line(0,-1){0.5}}
\put(5.4,-1.5){\line(1,0){11.5}}
\put(5.4,-1.5){\line(0,-1){0.5}}
\put(16.9,-1.5){\line(0,-1){0.5}}
\put(1.4,-3){\framebox(8.0,1)[c]
{\shortstack[c]{Inertial coordinates (Sec. 3.1)}}} 
\put(5.4,-3.0){\line(0,-1){0.5}}
\put(2.65,-3.5){\line(1,0){5.5}}
\put(2.65,-3.5){\line(0,-1){0.5}}
\put(8.15,-3.5){\line(0,-1){0.5}}
\put(0.0,-5){\framebox(5.3,1)[c]
{\shortstack[c]{Lagrangian (3.1.1)}}} 
\put(5.5,-5){\framebox(5.3,1)[c]
{\shortstack[c]{Hamiltonian (3.1.2)}}} 
\put(12.9,-3){\framebox(8,1)[c]
{\shortstack[c]{Relative coordinates (Sec. 3.2)}}} 
\put(16.9,-3.0){\line(0,-1){0.5}}
\put(14.15,-3.5){\line(1,0){5.5}}
\put(14.15,-3.5){\line(0,-1){0.5}}
\put(19.65,-3.5){\line(0,-1){0.5}}
\put(11.5,-5){\framebox(5.3,1)[c]
{\shortstack[c]{Lagrangian (3.2.1)}}} 
\put(17,-5){\framebox(5.3,1)[c]
{\shortstack[c]{Hamiltonian (3.2.2)}}} 
\put(6.15,-7.0){\framebox(10,1)[c]
{\shortstack[c]{Discrete time (Sec. 4)}}}
\put(11.15,-7){\line(0,-1){0.5}}
\put(5.4,-7.5){\line(1,0){11.5}}
\put(5.4,-7.5){\line(0,-1){0.5}}
\put(16.9,-7.5){\line(0,-1){0.5}}
\put(1.4,-9){\framebox(8,1)[c]
{\shortstack[c]{Inertial coordinates (Sec. 4.1)}}} 
\put(5.4,-9.0){\line(0,-1){0.5}}
\put(2.65,-9.5){\line(1,0){5.5}}
\put(2.65,-9.5){\line(0,-1){0.5}}
\put(8.15,-9.5){\line(0,-1){0.5}}
\put(0.0,-11){\framebox(5.3,1)[c]
{\shortstack[c]{Lagrangian (4.1.1)}}} 
\put(5.5,-11){\framebox(5.3,1)[c]
{\shortstack[c]{Hamiltonian (4.1.2)}}} 
\put(12.9,-9){\framebox(8,1)[c]
{\shortstack[c]{Relative coordinates (Sec. 4.2)}}} 
\put(16.9,-9.0){\line(0,-1){0.5}}
\put(14.15,-9.5){\line(1,0){5.5}}
\put(14.15,-9.5){\line(0,-1){0.5}}
\put(19.65,-9.5){\line(0,-1){0.5}}
\put(11.5,-11){\framebox(5.3,1)[c]
{\shortstack[c]{Lagrangian (4.2.1)}}} 
\put(17,-11){\framebox(5.3,1)[c]
{\shortstack[c]{Hamiltonian (4.2.2)}}} 
\end{picture}}
}
\caption{Eight types of equations of motion for the full body
problem}\label{fig:8eom}
\end{figure}
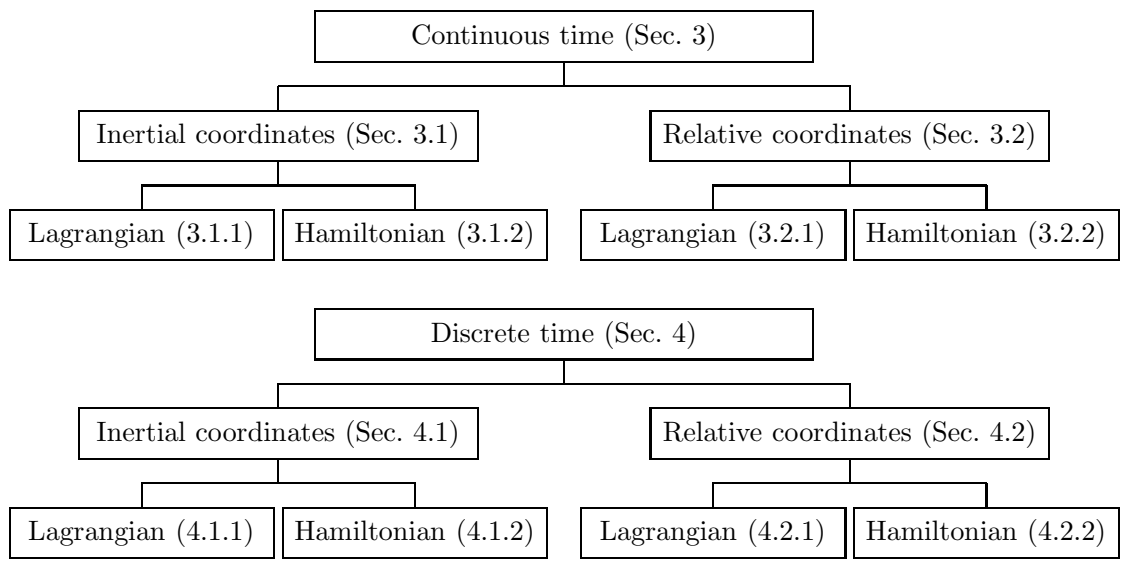

Continuous equations of motion for the full body problem are given
in~\cite{jo:macie} without any formal derivation of the equations.
We show, in this paper, that the equations of motion for the full
body problem can be derived from Hamilton's principle. The proper
form for the variations of Lie group elements in the configuration
space lead to a systematic derivation of the equations of motion.

This paper develops discrete variational equations of motion for the
full body model following a similar variational approach but carried
out within a discrete time framework. Since numerical integrators
are derived from the discrete Hamilton's principle, they exhibit
symplectic and momentum preserving properties, and good energy
behavior. They are also constructed to conserve the Lie group
geometry on the configuration space. Numerical simulation results
for the interaction of two rigid dumbbell models are given.

This paper is organized as follows. The basic idea of the
variational integrator is given in section 2. The continuous
equations of motion and variational integrators are developed in
section 3 and 4. Numerical simulation results are given in section
5. An appendix contains several technical details that are utilized
in the main development.

\section{Background}

\subsection{Hamilton's principle and variational integrators}
\begin{figure}[htbp]
\resizebox{1.00\linewidth}{!}{
\setlength{\unitlength}{1.65em}\centering\footnotesize
\begin{picture}(22.3,14)(0,-14)
\put(0.0,-2.0){\framebox(5.3,2.0)[c]
{\shortstack[c]{Configuration Space\\$(q,\dot{q})\in TQ$}}}
\put(2.65,-2.0){\vector(0,-1){1.0000}}
\put(0.0,-5.0){\framebox(5.3,2.0)[c]
{\shortstack[c]{Lagrangian\\$L(q,\dot{q})$}}}
\put(2.65,-5.0){\vector(0,-1){1.0000}}
\put(2.65,-5.5){\line(1,0){5.5}}
\put(8.15,-5.5){\vector(0,-1){3.5}}
\put(0.0,-8.0){\framebox(5.3,2.0)[c]
{\shortstack[c]{Action Integral\\$\mathfrak{G}=\int_{t_0}^{t_f}
L(q,\dot{q})\, dt$}}} \put(2.65,-8.0){\vector(0,-1){1.0000}}
\put(0.0,-11.0){\framebox(5.3,2.0)[c]
{\shortstack[c]{Variation\\$\delta\mathfrak{G}=\frac{d}{d\epsilon}\mathfrak{G}^\epsilon=0$}}}
\put(2.65,-11.0){\vector(0,-1){1.0000}}
\put(0.0,-14.0){\framebox(5.3,2.0)[c]
{\shortstack[c]{Euler-Lagrange Eqn.\\$\frac{d}{dt}\deriv{L}{\dot
q}-\deriv{L}{q}=0$}}}
\put(5.5,-11.0){\framebox(5.3,2.0)[c]
{\shortstack[c]{Legendre transform.\\$p=\mathbb{F}L(q,\dot{q})$}}}
\put(8.15,-11.0){\vector(0,-1){1.0000}}
\put(5.5,-14.0){\framebox(5.3,2.0)[c]
{\shortstack[c]{Hamilton's Eqn.\\$\dot{q}=H_p,\,\dot{p}=-H_q$}}}
\put(11.5,-2.0){\framebox(5.3,2.0)[c]
{\shortstack[c]{Configuration Space\\$(q_k, q_{k+1})\in Q\times
Q$}}} \put(14.15,-2.0){\vector(0,-1){1.0000}}
\put(11.5,-5.0){\framebox(5.3,2.0)[c]
{\shortstack[c]{Discrete Lagrangian\\$L_{d}(q_k,q_{k+1})$}}}
\put(14.15,-5.0){\vector(0,-1){1.0000}}
%
\put(14.15,-5.5){\line(1,0){5.5}}
\put(19.65,-5.5){\vector(0,-1){3.5}}
\put(11.5,-8.0){\framebox(5.3,2.0)[c]
{\shortstack[c]{Action Sum\\$\mathfrak{G}_d=\sum
L_{d}(q_k,q_{k+1})$}}} \put(14.15,-8.0){\vector(0,-1){1.0000}}
\put(11.5,-11.0){\framebox(5.3,2.0)[c]
{\shortstack[c]{Variation\\$\delta\mathfrak{G}_d=\frac{d}{d\epsilon}\mathfrak{G}_d^\epsilon=0$}}}
\put(14.15,-11.0){\vector(0,-1){1.0000}}
\put(11.5,-14.0){\framebox(5.3,2.0)[c]
{\shortstack[c]{Dis. E-L Eqn.
\\\scriptsize{$\mathbf{D}_{q_k}L_{d_{k-1}}\!+\!\mathbf{D}_{q_k}L_{d_k}=0$}}}}
\put(17,-11.0){\framebox(5.3,2.0)[c]
{\shortstack[c]{Legendre
transform.\\$p_k=\mathbb{F}L(q,\dot{q})\big|_k$}}}
\put(19.65,-11.0){\vector(0,-1){1.0000}}
\put(17,-14.0){\framebox(5.3,2.0)[c]
{\shortstack[c]{Dis. Hamilton's Eqn.
\\\scriptsize{$p_k=-\mathbf{D}_{q_k}L_{d_k},$}\\
\scriptsize{$p_{k+1}=\mathbf{D}_{q_{k+1}}L_{d_k}$}}}}
\end{picture}
}
\caption{Procedures to derive the continuous and discrete equations
of motion}\label{fig:el}
\end{figure}
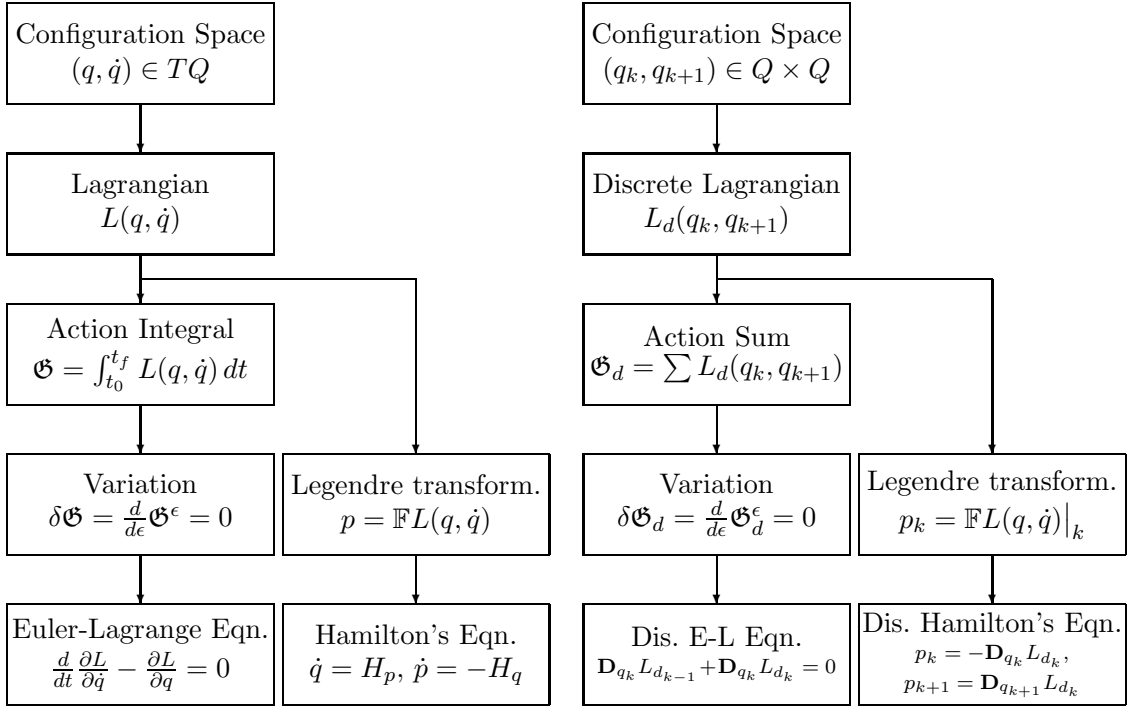

The procedure to derive the Euler-Lagrange equations and Hamilton's
equations from Hamilton's principle is shown in \reffig{el}.

When deriving the equations of motion, we first choose generalized
coordinates $q$, and a corresponding configuration space $Q$. We
then construct a Lagrangian from the kinetic and the potential
energy. An action integral $\mathfrak{G}=\int_{t_0}^{t_f} L(q,\dot
q) dt$ is defined as the path integral of the Lagrangian along a
time-parameterized trajectory. After taking the variation of the
action integral, and requiring it to be stationary, we obtain the
Euler-Lagrange equations. If we use the Legendre transformation
defined as
\begin{align}
p\cdot\delta\dot{q}&
=\mathbb{F}L(q,\dot q) \cdot\delta\dot{q},\nonumber\\
& = \frac{d}{d\epsilon} \bigg|_{\epsilon=0} L(q,\dot q+\epsilon
\delta\dot{q}),\label{eqn:FL}
\end{align}
where $\delta\dot{q}\in T_q Q$, then we obtain Hamilton's equations
in terms of momenta variables $p=\mathbb{F}L(q,\dot{q})\in T^*Q$.
These equations are equivalent to the Euler-Lagrange
equations~\cite{bk:marsden}.

There are two issues that arise in applying this procedure to the
full body problem. The first is that the configuration space for
each rigid body is the semi-direct product, $\SE=\Rset^3
\,\textcircled{s}\,\SO$, where $\SO$ is the Lie group of orthogonal
matrices with determinant 1. Therefore, variations should be
carefully chosen such that they respect the geometry of the
configuration space. For example, a varied rotation matrix
$R^\epsilon\in\SO$ can be expressed as
\begin{align*}
R^\epsilon = R e^{\epsilon\eta},
\end{align*}
where $\epsilon\in\Rset$, and $\eta\in\so$ is a variation
represented by a skew-symmetric matrix. The variation of the
rotation matrix $\delta R$ is the part of $R^\epsilon$ that is
linear in $\epsilon$:
\begin{align*}
R^\epsilon = R + \epsilon \delta R + \mathcal{O}(\epsilon^2).
\end{align*}
More precisely, $\delta R$ is given by
\begin{align}
\delta R =
\frac{d}{d\epsilon}\bigg|_{\epsilon=0}R^\epsilon=R\eta.\label{eqn:deltaR0}
\end{align}
The second issue is that reduced variables can be used to obtain
equations of motion expressed in relative coordinates. The
variations of the reduced variables are constrained as they must
arise from the variations of the unreduced variables while
respecting the geometry of the configuration space. The proper
variations of Lie group elements and reduced quantities are computed
while deriving the continuous equations of motion.

Generally, numerical integrators are obtained by approximating the
continuous Euler-Lagrange equation using a finite difference rule
such as $\dot{q}_k=(q_{k+1}-q_k)/h$, where $q_k$ denotes the value
of $q(t)$ at $t=hk$ for an integration step size $h\in\Rset$ and an
integer $k$. A variational integrator is derived by following a
procedure shown in the right column of \reffig{el}. When deriving a
variational integrator, the velocity phase space $(q,\dot q)\in TQ$
of the continuous Lagrangian is replaced by $(q_k,q_{k+1})\in
Q\times Q$, and the discrete Lagrangian $L_{d}$ is chosen such that
it approximates a segment of the action integral
\begin{align*}
L_d(q_k,q_{k+1})\approx\int_{0}^{h}
L\parenth{q_{k,k+1}(t),\dot{q}_{k,k+1}(t)} dt,
\end{align*}
where $q_{k,k+1}(t)$ is the solution of the Euler-Lagrange equation
satisfying boundary conditions $q_{k,k+1}(0)=q_k$ and
$q_{k,k+1}(h)=q_{k+1}$. Then, the discrete action sum
$\mathfrak{G}_d=\sum L_d(q_k, q_{k+1})$ approximates the action
integral $\mathfrak{G}$. Taking the variations of the action sum, we
obtain the discrete Euler Lagrange equation
\begin{align*}
\mathbf{D}_{q_k}L_d(q_{k-1},q_k)+\mathbf{D}_{q_k}
L_d(q_k,q_{k+1})=0,
\end{align*}
where $\mathbf{D}_{q_k}L_d$ denotes the partial derivative of $L_d$
with respect to $q_k$. This yields a discrete Lagrangian map
$F_{L_d}:(q_{k-1},q_k)\mapsto(q_k,q_{k+1})$. Using a discrete
analogue of the Legendre transformation, referred to as a discrete
fiber derivative $\mathbb{F}L_d:Q\times Q\rightarrow T^* Q$,
variational integrators can be expressed in Hamiltonian form as
\begin{align}
p_k & =-\mathbf{D}_{q_k} L_d (q_k,q_{k+1}),\label{eqn:dLtk}\\
p_{k+1} & = \mathbf{D}_{q_{k+1}} L_d (q_k,q_{k+1}).\label{eqn:dLtkp}
\end{align}
This yields a discrete Hamiltonian map
$\tilde{F}_{L_d}:(q_{k},p_k)\mapsto(q_{k+1},p_{k+1})$. A complete
development of variational integrators can be found
in~\cite{jo:marsden}.

\subsection{Notation}
Variables in the inertial and the body fixed coordinates are denoted
by lower-case and capital letters, respectively. A subscript $i$ is
used for variables related to the $i$th body. The relative variables
have no subscript and the $k$th discrete variables have the second
level subscript $k$. The letters $x,v,\Omega$ and $R$ are used to
denote the position, velocity, angular velocity and rotation matrix,
respectively.

The trace of $A\in\Rset^{n\times n}$ is denoted by
\begin{align*}
\tr{A}=\sum_{i=1}^{n}\, [A]_{ii},
\end{align*}
where $[A]_{ii}$ is the $i,i$th element of $A$. It can be shown that
\begin{gather}
\tr{AB}=\tr{BA}=\tr{B^TA^T}=\tr{A^TB^T},\label{eqn:trAB}\\
\tr{A^TB}=\sum_{p,q=1}^3 [A]_{pq} [B]_{pq},\label{eqn:trSum}\\
\tr{PQ}=0,\label{eqn:trskew}
\end{gather}
for matrices $A,B\in\Rset^{n\times n}$, a skew-symmetric matrix
$P=-P^T\in\Rset^{n\times n}$, and a symmetric matrix
$Q=Q^T\in\Rset^{n\times n}$.

The map $S:\Rset^3\mapsto\Rset^{3\times 3}$ is defined by the
condition that $S(x)y=x\times y$ for $x,y\in\Rset^3$. For
$x=(x_1,x_2,x_3)\in\Rset^3$, $S(x)$ is given by
\begin{align*}
S(x)=\begin{bmatrix}0&-x_3&x_2\\x_3&0&-x_1\\-x_2&x_1&0\end{bmatrix}.
\end{align*}
It can be shown that
\begin{gather}
S(x)^T=-S(x),\label{eqn:Sskew}\\
S(x\times y)=S(x)S(y)-S(y)S(x)=yx^T-xy^T,\label{eqn:Scross}\\
S(Rx)=RS(x)R^T,\label{eqn:SR}\\
S(x)^TS(x)=\parenth{x^T x} I_{3\times 3} - x x^T=\tr{xx^T}I_{3\times
3}-x x^T,\label{eqn:STS}
\end{gather}
for $x,y\in\Rset^3$ and $R\in\SO$.

Using homogeneous coordinates, we can represent an element of $\SE$
as follows:
\begin{align*}
\aSE{R}{x}\in\SE
\end{align*}
for $x\in\Rset^3$ and $R\in\SO$. Then, an action on $\SE$ is given
by the usual matrix multiplication in $\Rset^{4\times 4}$. For
example
\begin{align*}
\aSE{R_1}{x_1}\aSE{R_2}{x_2}=\aSE{R_1R_2}{R_1x_2+x_1}.
\end{align*}
for $x_1,x_2\in\Rset^3$ and $R_1,R_2\in\SO$.

The action of an element of $\SE$ on $\Rset^3$ can be expressed
using a matrix-vector product, once we identify $\Rset^3$ with
$\Rset^3\times\{1\}\subset\Rset^4$. In particular, we see from
\[ \aSE{R}{x_1} \begin{bmatrix} x_2 \\ 1\end{bmatrix} = \begin{bmatrix} Rx_2 + x_1\\ 1\end{bmatrix}\]
that $x_2\mapsto Rx_2+x_1$.

\section{Continuous time full body models}

In this section, the continuous time equations of motion for the
full body problem are derived in inertial and relative coordinates,
and are expressed in both Lagrangian and in Hamiltonian form.

We define $O-e_1e_2e_3$ as an inertial frame, and
$O_{\mathcal{B}_i}-E_{i_1}E_{i_2}E_{i_3}$ as a body fixed frame for
the $i$th body, $\mathcal{B}_i$. The origin of the $i$th body fixed
frame is located at the center of mass of body $\mathcal{B}_i$. The
configuration space of the $i$th rigid body is
$\SE=\Rset^3\,\textcircled{s}\,\SO$. We denote the position of the
center of mass of $\mathcal{B}_i$ in the inertial frame by
$x_i\in\Rset^3$,\label{no:xi} and we denote the attitude of
$\mathcal{B}_i$ by $R_i\in\SO$,\label{no:Ri} which is a rotation
matrix from the $i$th body fixed frame to the inertial frame.

\subsection{Inertial coordinates}

\textit{Lagrangian:} To derive the equations of motion, we first
construct a Lagrangian for the full body problem. Given
$(x_i,R_i)\in\SE$, the inertial position of a mass element of
$\mathcal{B}_i$ is given by $x_i+R_i\rho_i$, where
$\rho_i\in\Rset^3$ denotes the position of the mass element in the
body fixed frame. Then, the kinetic energy of $\mathcal{B}_i$ can be
written as
\begin{align*}
T_i & = \frac{1}{2}\int_{\mathcal{B}_i}
\|\dot{x}_i+\dot{R}_i\rho_i\|^2\, dm_i.
\end{align*}
Using the fact that $\int_{\mathcal{B}_i} \rho_i dm_i=0$ and the
kinematic equation $\dot{R}_i=R_iS(\Omega_i)$, the kinetic energy
$T_i$ can be rewritten as
\begin{align}
T_i(\dot{x}_i,\Omega_i)
 & =
\frac{1}{2}\int_{\mathcal{B}_i} \norm{\dot{x}_i}^2 +
\norm{S(\Omega_i)\rho_i}^2 dm_i,\nonumber\\
& = \frac{1}{2}m_i \norm{\dot{x}_i}^2 +
\frac{1}{2}\tr{S(\Omega_i)J_{d_i}S(\Omega_i)^T},\label{eqn:Ti}
\end{align}
where $m_i\in\Rset$ is the total mass of $\mathcal{B}_i$\label{no:mi},
$\Omega_i\in\Rset^3$ is the angular velocity of $\mathcal{B}_i$ in
the body fixed frame,\label{no:Omegai} and
$J_{d_i}=\int_{\mathcal{B}_i}\rho_i\rho_i^Tdm_i\in\Rset^{3\times
3}$\label{no:Jdi} is a nonstandard moment of inertia matrix. Using
\refeqn{STS}, it can be shown that the standard moment of inertia
matrix
$J_{i}=\int_{\mathcal{B}_i}S(\rho_i)^TS(\rho_i)dm_i\in\Rset^{3\times
3}$\label{no:Ji} is related to $J_{d_i}$ by
\begin{align*}
J_i=\tr{J_{d_i}}I_{3\times 3} - J_{d_i},
\end{align*}
and that the following condition holds for any $\Omega_i\in\Rset^3$
\begin{align}
S(J_i\Omega_i)=S(\Omega_i)J_{d_i}+J_{d_i}S(\Omega_i).\label{eqn:JdJ}
\end{align}
We first derive equations using the nonstandard moment of inertia
matrix $J_{d_i}$, and then express them in terms of the standard
moment of inertia $J_i$ by using \refeqn{JdJ}.

The gravitational potential energy $U:\SE^n\mapsto\Rset$ is given by
\begin{align}
U(x_1,x_2,\cdots,x_n,R_1,R_2,\cdots,R_n) =-\frac{1}{2}
\sum_{\substack{i,j=1\\i\neq
j}}^{n}\int_{\mathcal{B}_i}\int_{\mathcal{B}_j}
\frac{G dm_i dm_j}{\norm{x_i+R_i\rho_i - x_j -
R_j\rho_j}},\label{eqn:U}
\end{align}
where $G$ is the universal gravitational constant.

Then, the Lagrangian for $n$ full bodies can be written as
\begin{align}
L(\mathbf{x},\mathbf{\dot
x},\mathbf{R},\mathbf{\Omega})=\sum_{i=1}^{n} \frac{1}{2}m_i
\norm{\dot{x}_i}^2 + \frac{1}{2}\tr{S(\Omega_i)J_{d_i}S(\Omega_i)^T}
- U(\mathbf{x},\mathbf{R}),\label{eqn:L}
\end{align}
where bold type letters denote ordered $n$-tuples of variables. For
example, $\mathbf{x}\in(\Rset^3)^n$, $\mathbf{R}\in\SO^n$, and
$\mathbf{\Omega}\in(\Rset^3)^n$ are defined as
$\mathbf{x}=(x_1,x_2,\cdots,x_n)$,\label{no:bfx}
$\mathbf{R}=(R_1,R_2,\cdots,R_n)$,\label{no:bfR} and
$\mathbf{\Omega}=(\Omega_1,\Omega_2,\cdots,\Omega_n)$,
respectively.\label{no:bfOmega}

\textit{Variations of variables:} Since the configuration space is
$\SE^n$, the variations should be carefully chosen such that they
respect the geometry of the configuration space. The variations of
$x_i:[t_0,t_f]\mapsto\Rset^3$ and
$\dot{x}_i:[t_0,t_f]\mapsto\Rset^3$ are trivial, namely
\begin{align*}
x_i^\epsilon = x_i + \epsilon\delta x_i +
\mathcal{O}(\epsilon^2),\\
\dot{x}_i^\epsilon = \dot{x}_i + \epsilon\delta \dot{x}_i +
\mathcal{O}(\epsilon^2),
\end{align*}
where $\delta x_i:[t_0,t_f]\mapsto\Rset^3$, $\delta
\dot{x}_i:[t_0,t_f]\mapsto\Rset^3$ vanish at the initial time $t_0$
and at the final time $t_f$. The variation of
$R_i:[t_0,t_f]\mapsto\SO$, as given in \refeqn{deltaR0}, is
\begin{align}
\delta R_i = R_i \eta_i\label{eqn:delRi},
\end{align}
where $\eta_i:[t_0,t_f]\mapsto\so$ is a variation with values
represented by a skew-symmetric matrix ($\eta_i^T=-\eta_i$)
vanishing at $t_0$ and $t_f$. The variation of $\Omega_i$ can be
computed from the kinematic equation $\dot{R}_i=R_iS(\Omega_i)$ and
\refeqn{delRi} to be
\begin{align}
S(\delta \Omega_i) & = \frac{d}{d\epsilon}\bigg|_{\epsilon=0}
R_i^{\epsilon T}\dot{R}_i^\epsilon=
\delta R_i^T \dot{R}_i + R_i^T \delta \dot{R}_i,\nonumber\\
 & = -\eta_iS(\Omega_i) +
S(\Omega_i)\eta_i +\dot{\eta}_i .\label{eqn:delOmegai}
\end{align}

\subsubsection{Equations of motion: Lagrangian form}
If we take variations of the Lagrangian using \refeqn{delRi} and
\refeqn{delOmegai}, we obtain the equations of motion from
Hamilton's principle. We first take the variation of the kinetic
energy of $\mathcal{B}_i$.
\begin{align*}
\delta T_i & =
\frac{d}{d\epsilon}\bigg|_{\epsilon=0}T_i(\dot{x}_i+\epsilon\delta\dot{x}_i,\Omega_i+\epsilon\delta\Omega_i),\nonumber\\
& = m_i \dot{x}_i^T \delta\dot{x}_i +
\frac{1}{2}\tr{S(\delta\Omega_i)J_{d_i}S(\Omega_i)^T+S(\Omega_i)J_{d_i}S(\delta\Omega_i)^T}.
\end{align*}
Substituting \refeqn{delOmegai} into the above equation and using
\refeqn{trAB}, we obtain
\begin{align*}
\delta T_i & = m_i \dot{x}_i^T \delta\dot{x}_i +
\frac{1}{2}\tr{-\dot\eta_i
\braces{J_{d_i}S(\Omega_i)+S(\Omega_i)J_{d_i}}}\nonumber\\
& \quad + \frac{1}{2}\tr{
\eta_i\braces{S(\Omega_i)\braces{J_{d_i}S(\Omega_i)+S(\Omega_i)J_{d_i}}
-\braces{J_{d_i}S(\Omega_i)+S(\Omega_i)J_{d_i}}S(\Omega_i)}}.
\end{align*}
Using \refeqn{Scross} and \refeqn{JdJ}, $\delta T_i$ is given by
\begin{align}
\delta T_i & = m_i \dot{x}_i^T \delta\dot{x}_i
+\frac{1}{2}\tr{-\dot{\eta}_iS(J_i\Omega_i)+\eta_iS(\Omega_i\times
J_i\Omega_i)}.\label{eqn:delTi}
\end{align}

The variation of the potential energy is given by
\begin{align*}
\delta U & = \frac{d}{d\epsilon}\bigg|_{\epsilon=0}
U(\mathbf{x}+\epsilon\delta\mathbf{x},\mathbf{R}+\epsilon\delta\mathbf{R}),
\end{align*}
where $\delta\mathbf{x}=(\delta x_1,\delta x_2,\cdots,\delta x_n)$,
$\delta \mathbf{R}=(\delta R_1,\delta R_2,\cdots,\delta R_n)$. Then,
$\delta U$ can be written as
\begin{align}
\delta U &
=\sum_{i=1}^n\parenth{\sum_{p=1}^3\deriv{U}{[x_i]_p}[\delta x_i]_p
+ \sum_{p,q=1}^{3}\deriv{U}{[R_i]_{pq}}[R_i\eta_i]_{pq}},\nonumber\\
& =\sum_{i=1}^n\parenth{\deriv{U}{x_i}^T\delta
x_i-\tr{\eta_iR^T_i\deriv{U}{R_i}}},\label{eqn:delU}
\end{align}
where $[A]_{pq}$ denotes the $p,q$th element of a matrix $A$, and
$\deriv{U}{x_i},\,\deriv{U}{R_i}$ are given by
$[\deriv{U}{x_i}]_p=\deriv{U}{[x_i]_p}$,
$[\deriv{U}{R_i}]_{pq}=\deriv{U}{[R_i]_{pq}}$, respectively. The
variation of the Lagrangian has the form \begin{align} \delta L =
\sum_{i=1}^{n} \delta T_i - \delta U,\label{eqn:delL}
\end{align}
which can be written more explicitly by using \refeqn{delTi} and
\refeqn{delU}.

The action integral is defined to be
\begin{align}
\mathfrak{G} = \int_{t_0}^{t_f} L(\mathbf{x},\mathbf{\dot
x},\mathbf{R},\mathbf{\Omega})\,dt.\label{eqn:caction}
\end{align}
The variation of the action integral can be written as
\begin{align*}
\delta \mathfrak{G} = \sum_{i=1}^{n} \int_{t_0}^{t_f} & 
m_i \dot{x}_i^T \delta\dot{x}_i-\deriv{U}{x_i}^T\delta x_i
+\frac{1}{2}\tr{-\dot{\eta}_iS(J_i\Omega_i)+
\eta_i\braces{S(\Omega_i\times J_i\Omega_i)
+2R^T_i\deriv{U}{R_i}}}\,dt.
\end{align*}
Using integration by parts,
\begin{align*}
\delta \mathfrak{G} & = \sum_{i=1}^{n}\; 
m_i\dot{x}_i^T \delta x_i \bigg|_{t_0}^{t_f}
-\frac{1}{2}\tr{\eta_iS(J_i\Omega_i)}\bigg|_{t_0}^{t_f}
+\int_{t_0}^{t_f} -m_i\ddot{x}_i^T \delta x_i
+\frac{1}{2}\tr{\eta_iS(J_i\dot\Omega_i)}\,dt\\
& \quad +  \sum_{i=1}^{n} \int_{t_0}^{t_f}
-\deriv{U}{x_i}^T\delta x_i
+\frac{1}{2}\tr{ \eta_i\braces{S(\Omega_i\times J_i\Omega_i)
+2R^T_i\deriv{U}{R_i}}}\,dt.
\end{align*}
Using the fact that $\delta x_i$ and $\eta_i$ vanish at $t_0$ and
$t_f$, the first two terms of the above equation vanish. Then,
$\delta\mathfrak{G}$ is given by
\begin{align*}
\delta \mathfrak{G} & = \sum_{i=1}^{n} \int_{t_0}^{t_f} 
- \delta x_i^T \braces{m_i\ddot{x}_i
+\deriv{U}{x_i}}+\frac{1}{2}\tr{\eta_i\braces{S(J_i\dot{\Omega}_i+\Omega_i\times
J_i\Omega_i)+2R^T_i\deriv{U}{R_i}}} \,dt.
\end{align*}

From Hamilton's principle, $\delta\mathfrak{G}$ should be zero for
all possible variations $\delta x_i:[t_0,t_f]\mapsto\Rset^3$ and
$\eta_i:[t_0,t_f]\mapsto\so$. Therefore, the expression in the first
brace should be zero. Furthermore, since $\eta_i$ is skew symmetric,
we have by \refeqn{trskew}, that  the expression in the second brace
should be symmetric. Then, we obtain the continuous equations of
motion as
\begin{gather}
m_i \ddot{x}_i = -\deriv{U}{x_i},\nonumber\\
S(J_i\dot{\Omega}_i+\Omega_i\times
J_i\Omega_i)+2R^T_i\deriv{U}{R_i}=
S(J_i\dot{\Omega}_i+\Omega_i\times
J_i\Omega_i)^T+2\deriv{U}{R_i}^TR.\label{eqn:SS}
\end{gather}
Using \refeqn{Sskew}, we can simplify \refeqn{SS} to be
\begin{gather*}
S(J_i\dot{\Omega}_i+\Omega_i\times J_i\Omega_i)
=\deriv{U}{R_i}^TR-R^T_i\deriv{U}{R_i}.
\end{gather*}
Note that the right hand side expression in the above equation is
also skew symmetric. Then, the moment due to the gravitational
potential on the $i$th body, $M_i\in\Rset^3$ is given by
$S(M_i)=\deriv{U}{R_i}^TR_i-R^T_i\deriv{U}{R_i}$. This moment $M_i$
can be expressed more explicitly as the following computation shows.
\begin{align*}
S(M_i) & = \deriv{U}{R_i}^TR_i-R^T_i\deriv{U}{R_i},\\
 & =
\begin{bmatrix} & &
\\u^T_{ri_1}&u^T_{ri_2}&u^T_{ri_3}\\& & \end{bmatrix}
\begin{bmatrix}&r_{i_1}&\\&r_{i_2}&\\&r_{i_3}\end{bmatrix}
-\begin{bmatrix} & &
\\r_{i_1}^T&r_{i_2}^T&r_{i_3}^T\\& & \end{bmatrix}
\begin{bmatrix}&u_{ri_1}&\\&u_{ri_2}&\\&u_{ri_3}\end{bmatrix},\\
& = \parenth{u^T_{ri_1}r_{i_1}-r_{i_1}^Tu_{ri_1}}
+\parenth{u^T_{ri_2}r_{i_2}-r_{i_2}^Tu_{ri_2}}
+\parenth{u^T_{ri_3}r_{i_3}-r_{i_3}^Tu_{ri_3}},
\end{align*}
where $r_{i_p},u_{ri_p}\in\Rset^{1\times 3}$ are the $p$th row
vectors of $R_i$ and $\deriv{U}{R_i}$, respectively. Substituting
$x=r_{i_p}^T$, $y=u_{ri_p}^T$ into \refeqn{Scross}, we obtain
\begin{align}
S(M_i) & = S(r_{i_1}\times u_{ri_1})+S(r_{i_2}\times
u_{ri_2})+S(r_{i_3}\times u_{ri_3}),\nonumber\\
&=S(r_{i_1}\times u_{ri_1}+r_{i_2}\times u_{ri_2}+r_{i_3}\times
u_{ri_3}),\label{eqn:SMi}
\end{align}
Then, the gravitational moment $M_i$ is given by
\begin{align}
M_i=r_{i_1}\times u_{ri_1}+r_{i_2}\times u_{ri_2}+r_{i_3}\times
u_{ri_3}.\label{eqn:Mi}
\end{align}\label{no:Mi}

In summary, \textit{the continuous equations of motion for the full
body problem, in Lagrangian form,} can be written for
$i\in(1,2,\cdots,n)$ as
\begin{gather}
\dot{v}_i=-\frac{1}{m_i}\deriv{U}{x_i},\label{eqn:vidot}\\
J_i\dot{\Omega}_i+\Omega_i\times J_i\Omega_i = M_i,\\
\dot{x}_i=v_i,\\
\dot{R}_i=R_iS(\Omega_i),\label{eqn:Ridot}
\end{gather}
where the translational velocity $v_i\in\Rset^3$ is defined as
$v_i=\dot{x}_i$.\label{no:vi}

\subsubsection{Equations of motion: Hamiltonian form}
We denote the linear and angular momentum of the $i$th body
$\mathcal{B}_i$ by $\gamma_i\in\Rset^3$\label{no:gammai}, and
$\Pi_i\in\Rset^3$\label{no:Pii}, respectively. They are related to
the linear and angular velocities by $\gamma_i = m_i v_i$, and
$\Pi_i = J_i \Omega_i$. Then, the equations of motion can be
rewritten in terms of the momenta variables. \textit{The continuous
equations of motion for the full body problem, in Hamiltonian form,}
can be written for $i\in(1,2,\cdots,n)$ as
\begin{gather}
\dot{\gamma}_i=-\deriv{U}{x_i},\label{eqn:gamidot}\\
\dot{\Pi}_i+\Omega_i\times\Pi_i = M_i,\\
\dot{x}_i=\frac{\gamma_i}{m_i},\\
\dot{R}_i=R_iS(\Omega_i).\label{eqn:Ridoth}
\end{gather}

\subsection{Relative coordinates}
The motion of the full rigid bodies depends only on the relative
positions and the relative attitudes of the bodies. This is a
consequence of the property that the gravitational potential can be
expressed using only these relative variables. Physically, this is
related to the fact that the total linear momentum and the total
angular momentum about the mass center of the bodies are conserved.
Mathematically, the Lagrangian is invariant under a left action of
an element of $\SE$. So, it is natural to express the equations of
motion in one of the body fixed frames. In this section, the
equations of motion for the full two body problem are derived in
relative coordinates. This result can be readily generalized to the
$n$ body problem.

\textit{Reduction of variables:} In~\cite{jo:macie}, the reduction
is carried out in stages, by first reducing position variables in
$\Rset^3$, and then reducing attitude variables in $\SO$. This is
equivalent to directly reducing the position and the attitude
variables in $\SE$ in a single step, which is a result that can be
explained by the general theory of Lagrangian reduction by
stages~\cite{jo:cemara}. The reduced position and the reduced
attitude variables are the relative position and the relative
attitude of the first body with respect to the second body. In other
words, the variables are reduced by applying the inverse of
$(R_2,x_2)\in\SE$ given by $(R_2^T,-R_2^Tx_2)\in\SE$, to the
following homogeneous transformations:
\begin{align}
\aSE{R_2^T}{-R_2^Tx_2}\parenth{\aSE{R_1}{x_1},\aSE{R_2}{x_2}}& =
\parenth{\aSE{R_2^TR_1}{R_2^T(x_1-x_2)},\aSE{R_2^TR_2}{R_2^T(x_2-x_2)}},\nonumber\\
& =
\parenth{\aSE{R_2^TR_1}{R_2^T(x_1-x_2)},\aSE{I_{3\times 3}}{0}}.\label{eqn:HT}
\end{align}
This motivates the definition of the reduced variables as
\begin{align}
X & = R_2^T(x_1-x_2),\label{eqn:X}\\
R & = R_2^TR_1,\label{eqn:R}
\end{align}
where $X\in\Rset^3$ is the relative position of the first body with
respect to the second body expressed in the second body fixed frame,
and $R\in\SO$ is the relative attitude of the first body with
respect to the second body. The corresponding linear and angular
velocities are also defined as
\begin{align}
V & = R_2^T(\dot{x}_1-\dot{x}_2),\\
\Omega & = R\Omega_1,
\end{align}
where $V\in\Rset^3$ represents the relative velocity of the first
body with respect to the second body in the second body fixed frame\label{no:V},
and $\Omega\in\Rset^3$ is the angular velocity of the first body
expressed in the second body fixed frame\label{no:Omega}. Here, the capital letters
denote variables expressed in the second body fixed frame.

For convenience, we denote the inertial position and the inertial
velocity of the second body, expressed in the second body fixed
frame by $X_2,V_2\in\Rset^3$:
\begin{align}
X_2 & = R_2^T x_2,\label{eqn:X2}\\
V_2 & = R_2^T \dot{x}_2.\label{eqn:V2}
\end{align}

\textit{Reduced Lagrangian:} The equations of motion in relative
coordinates are derived in the same way used to derive the equations
in the inertial frame. We first construct a reduced Lagrangian. The
reduced Lagrangian $l$ is obtained by expressing the original
Lagrangian \refeqn{L} in terms of the reduced variables. The kinetic
energy is given by
\begin{align*}
T_1 + T_2 & =
\frac{1}{2}m_1\norm{\dot{x}_1}^2 +\frac{1}{2}m_2\norm{\dot{x}_2}^2
+\frac{1}{2}\tr{S(\Omega_1)J_{d_1}S(\Omega_1)^T}
+\frac{1}{2}\tr{S(\Omega_2)J_{d_2}S(\Omega_2)^T},\\
& =\frac{1}{2}m_1\norm{V+V_2}^2 +\frac{1}{2}m_2\norm{V_2}^2
+\frac{1}{2}\tr{S(\Omega)J_{d_R}S(\Omega)^T}
+\frac{1}{2}\tr{S(\Omega_2)J_{d_2}S(\Omega_2)^T},
\end{align*}
where \refeqn{SR} is used, and $J_{d_{R}}\in\Rset^{3\times 3}$ is
defined as $J_{d_{R}}=RJ_{d_{1}}R^T$, which is an expression of the
nonstandard moment of inertia matrix of the first body
with respect to the second body fixed frame.\label{no:JdR} Note that $J_{d_{R}}$
is not a constant matrix. Using \refeqn{SR}, it can be shown that
$J_{d_{R}}$ also satisfies a property similar to \refeqn{JdJ},
namely
\begin{align}
S(J_{R}\Omega)=S(\Omega)J_{d_{R}}+J_{d_{R}}S(\Omega),\label{eqn:JdRJR}
\end{align}
where $J_R=RJ_1R^T\in\Rset^{3\times 3}$ is the standard moment of
inertia matrix of the first body with respect to the second body
fixed frame.\label{no:JR}

Using \refeqn{U}, the gravitational potential $U$ of the full two
bodies is given by
\begin{align*}
U(x_1,x_2,R_1,R_2) =- \int_{\mathcal{B}_1}\int_{\mathcal{B}_2}
\frac{G dm_1 dm_2}{\norm{x_1+R_1\rho_1 - x_2 - R_2\rho_2}},
\end{align*}
and it is invariant under an action of an element of $\SE$.
Therefore, the gravitational potential can be written as a function
of the relative variables only. By applying the inverse of
$(R_2,x_2)\in\SE$ as given in \refeqn{HT}, we obtain
\begin{align*}
U(x_1,x_2,R_1,R_2) & = U(R_2^T(x_1-x_2),0,R_2^TR_1,I_{3\times 3}),\\
& = - \int_{\mathcal{B}_1}\int_{\mathcal{B}_2}
\frac{G dm_1 dm_2}{\norm{R_2^T(x_1-x_2)+R_2^T R_1\rho_1 - I_{3\times
3}
\rho_2}},\\
& = - \int_{\mathcal{B}_1}\int_{\mathcal{B}_2}
\frac{G dm_1 dm_2}{\norm{X+R\rho_1 - \rho_2}},
\\ & \triangleq U(X,R).
\end{align*}
Here we abuse notation slightly by using the same letter $U$ to
denote the gravitational potential as a function of the relative
variables.

Then, the reduced Lagrangian $l$ is given by
\begin{align}
l (R,X,\Omega,V,\Omega_2,V_2)& = \frac{1}{2} m_1 \|V+V_2\|^2 +
\frac{1}{2} m_2 \|V_2\|^2 \nonumber
\\ & \quad+
\frac{1}{2}\tr{S(\Omega)J_{d_{R}}S(\Omega)^T} + \frac{1}{2}
\tr{S(\Omega_2)J_{d_2} S(\Omega_2)^T}- U(X,R).\label{eqn:l}
\end{align}

\textit{Variations of reduced variables:} The variations of the
reduced variables must be restricted to those that can arise from
the variations of the original variables. For example, the variation
of the relative attitude $R$ is given by
\begin{align*}
\delta R & = \frac{d}{d\epsilon}\bigg|_{\epsilon=0}
R_2^{\epsilon T}R_1^\epsilon,\\
& = \delta R_2^TR_1 + R_2^T \delta R_1.
\end{align*}
Substituting \refeqn{delRi} into the above equation,
\begin{align*}
\delta R & = -\eta_2R_2^TR_1 + R_2^T R_1\eta_1,\\
& = -\eta_2R + \eta R,
\end{align*}
where a reduced variation $\eta\in\so$ is defined as
$\eta=R\eta_1R^T$. The variations of other reduced variables can be
obtained in a similar way. The detailed derivations are given in
\ref{apprv}, and we summarize the results as follows:
\begin{align}
\delta R & = \eta R - \eta_2 R,\label{eqn:delR}\\
\delta X & = \chi - \eta_2 X,\label{eqn:delX}\\
S(\delta\Omega) & = \dot\eta -S(\Omega)\eta+\eta
S(\Omega)+S(\Omega)\eta_2-\eta_2
S(\Omega)+S(\Omega_2)\eta-\eta S(\Omega_2),\label{eqn:delOmega}\\
\delta V & = \dot \chi + S(\Omega_2) \chi -\eta_2 V,\label{eqn:delV}\\
S(\delta\Omega_2) & = \dot\eta_2 +S(\Omega_2)\eta_2-\eta_2 S(\Omega_2),\label{eqn:delOmega2}\\
\delta V_2 & = \dot \chi_2 + S(\Omega_2) \chi_2 -\eta_2
V_2,\label{eqn:delV2}
\end{align}
where $\chi,\chi_2 :[t_0,t_f]\mapsto\Rset^3$ and $\eta,\eta_2
:[t_0,t_f]\mapsto\so$ are variations that vanish at the end points.
These Lie group variations are the key elements required to obtain
the equations of motion in relative coordinates.

\subsubsection{Equations of motion: Lagrangian form}
The reduced equations of motion can be computed from the reduced
Lagrangian using  the reduced Hamilton's principle. By taking the
variation of the reduced Lagrangian \refeqn{l} using the constrained
variations given by \refeqn{delR} through \refeqn{delV2}, we can
obtain the equations of motion in the relative coordinates.

Following a similar process to the derivation of $\delta T_i, \delta
U$ as in \refeqn{delTi} and \refeqn{delU}, the variation of the
reduced Lagrangian $\delta l$ can be obtained as
\begin{align}
\delta l & = \dot\chi^T \bracket{m_1(V+V_2)} - \chi^T\bracket{m_1
\Omega_2 \times (V+V_2)}\nonumber\\
& \quad + \dot\chi_2^T
\bracket{m_1(V+V_2)+m_2V_2}-\chi_2^T\bracket{m_1 \Omega_2 \times (V+V_2)+m_2 \Omega_2 \times V_2}\nonumber\\
& \quad +\frac{1}{2}\tr{-\dot\eta S(J_R \Omega)+\eta
S(\Omega_2\times J_R
\Omega)}+\frac{1}{2}\tr{-\dot\eta_2S(J_2\Omega_2)+\eta_2
S(\Omega_2\times J_2\Omega_2)}\nonumber\\
&\quad - \chi^T \deriv{U}{X} + \tr{\eta_2X\deriv{U}{X}^T} +
\tr{\eta_2R\deriv{U}{R}^T-\eta R \deriv{U}{R}^T},\label{eqn:dell}
\end{align}
where we used the identities \refeqn{Scross}, \refeqn{JdJ} and
\refeqn{JdRJR}, and the constrained variations \refeqn{delR} through
\refeqn{delV2}.

The action integral in terms of the reduced Lagrangian is
\begin{align}
\mathfrak{G} = \int_{t_0}^{t_f} l (R,X,\Omega,V,\Omega_2,V_2)
\,dt.\label{eqn:craction}
\end{align}
Using integration by parts together with the fact that
$\chi,\chi_2,\eta$ and $\eta_2$ vanish at $t_0$ and $t_f$, the
variation of the action integral can be expressed from \refeqn{dell}
as
\begin{align*}
\delta\mathfrak{G} & = - \int_{t_0}^{t_f} \chi^T \braces{m_1(\dot
V+\dot V_2)+m_1
\Omega_2 \times (V+V_2)+\deriv{U}{X}}\,dt\\
& \quad -\int_{t_0}^{t_f} \chi_2^T
\braces{m_1(\dot V+\dot V_2)+m_2\dot V_2+m_1 \Omega_2 \times (V+V_2)+m_2 \Omega_2 \times V_2}\,dt\\
& \quad +\frac{1}{2}\int_{t_0}^{t_f}\tr{\eta \braces{S(\dot{(J_R\Omega)}+\Omega_2\times J_R\Omega)-2R \deriv{U}{R}^T}}\,dt\\
& \quad
+\frac{1}{2}\int_{t_0}^{t_f}\tr{\eta_2\braces{S(J_2\dot{\Omega}_2+\Omega_2\times
J_2\Omega_2) +2X\deriv{U}{X}^T+2R\deriv{U}{R}^T}}\,dt.
\end{align*}

From the reduced Hamilton's principle, $\delta\mathfrak{G}=0$ for
all possible variations $\chi,\chi_2:[t_0,t_f]\mapsto\Rset^3$ and
$\eta,\eta_2:[t_0,t_f]\mapsto\so$ that vanish at $t_0$ and $t_f$.
Therefore, in the above equation, the expressions in the first two
braces should be zero and the expressions in the last two braces
should be symmetric since $\eta,\eta_2$ are skew symmetric. Then, we
obtain the following equations of motion,
\begin{gather}
m_1(\dot V+\dot V_2)+m_1 \Omega_2 \times (V+V_2)=-\deriv{U}{X},\label{eqn:Vdot0}\\
m_2\dot V_2+m_2\Omega_2 \times V_2=\deriv{U}{X},\label{eqn:V2dot0}\\
S(\dot{(J_R\Omega)}+\Omega_2\times J_R\Omega)=-S(M),\nonumber\\
S(J_2\dot{\Omega}_2+\Omega_2\times
J_2\Omega_2)=\deriv{U}{X}X^T-X\deriv{U}{X}^T+S(M),\label{eqn:Pi2dot0}
\end{gather}
where $M\in\Rset^3$ is defined by the relation
$S(M)=\deriv{U}{R}R^T-R \deriv{U}{R}^T$. By a procedure analogous to
the derivation of \refeqn{SMi}, $M$ can be written as
\begin{align}
M = r_1 \times u_{r_1} + r_2 \times u_{r_2} + r_3 \times
u_{r_3},\label{eqn:M}
\end{align}
where $r_p,u_{r_p}\in\Rset^{3}$ are the $p$th column vectors of $R$
and $\deriv{U}{R}$, respectively.

Equation \refeqn{Vdot0} can be simplified using \refeqn{V2dot0} as
\begin{align*}
\dot{V}+\Omega_2\times V = -\frac{m_1+m_2}{m_1m_2}\deriv{U}{X}.
\end{align*}
For reconstruction of the motion of the second body, it is natural
to express the motion of the second body in the inertial frame.
Since $\dot{V}_2=\dot{R}_2^T\dot{x}_2+R_2^T\ddot{x}_2
=-S(\Omega_2)V+R_2^T\dot{v}_2$, \refeqn{V2dot0} can be written as
\begin{align*}
m_2 R_2^T\dot{v}_2 = \deriv{U}{X}.
\end{align*}
Equation \refeqn{Pi2dot0} can be simplified using the property
$\deriv{U}{X}X^T-X\deriv{U}{X}^T=S(X\times\deriv{U}{X})$ from
\refeqn{Scross}. The kinematics equations for $\dot{R}$ and
$\dot{X}$ can be derived in a similar way.

In summary, \textit{the continuous equations of relative motion for
the full two body problem, in Lagrangian form,} can be written as
\begin{gather}
\dot{V}+\Omega_2\times V = -\frac{1}{m}\deriv{U}{X},\label{eqn:Vdot}\\
\dot{(J_R\Omega)}+\Omega_2\times J_R\Omega=-M,\\
J_2\dot{\Omega}_2+\Omega_2\times J_2\Omega_2=X \times
\deriv{U}{X}+M,\label{eqn:Pi2dot}\\
\dot{X}+\Omega_2 \times X=V,\\
\dot{R}=S(\Omega)R-S(\Omega_2)R,\label{eqn:Rdot}
\end{gather}
where $m=\frac{m_1m_2}{m_1+m_2}$.\label{no:m} The following equations can be
used for reconstruction of the motion of the second body in the
inertial frame:
\begin{gather}
\dot{v}_2 = \frac{1}{m_2}R_2 \deriv{U}{X},\label{eqn:v2dot}\\
\dot{x}_2 = v_2,\\
\dot{R_2}=R_2 S(\Omega_2).\label{eqn:R2dot}
\end{gather}

These equations are equivalent to those given in \cite{jo:macie}.
However, \refeqn{v2dot} is not given in \cite{jo:macie}. Equations
\refeqn{Vdot} though \refeqn{R2dot} give a complete set of equations
for the reduced dynamics and reconstruction. Furthermore, they are
derived systematically in the context of geometric mechanics using
proper variational formulas given in \refeqn{delR} through
\refeqn{delV2}. This result can be readily generalized for $n$
bodies.

\subsubsection{Equations of motion: Hamiltonian form}
Define the linear momenta $\Gamma,\gamma_2\in\Rset^3$,\label{no:Gamma} and the
angular momenta $\Pi,\Pi_2\in\Rset^3$\label{no:Pi} as
\begin{align*}
\Gamma & = m V,\\
\gamma_2 & = m v_2,\\
\Pi& =J_R\Omega=RJ_1\Omega_1,\\
\Pi_2&=J_2\Omega_2.
\end{align*}
Then, the equations of motion can be rewritten in terms of these
momenta variables. \textit{The continuous equations of relative
motion for the full two body problem, in Hamiltonian form,} can be
written as
\begin{gather}
\dot{\Gamma}+\Omega_2\times \Gamma = -\deriv{U}{X},\label{eqn:Gamdot}\\
\dot{\Pi}+\Omega_2\times \Pi=-M,\\
\dot{\Pi}_2+\Omega_2\times \Pi_2=X \times
\deriv{U}{X}+M,\\
\dot{X}+\Omega_2 \times X=\frac{\Gamma}{m},\\
\dot{R}=S(\Omega)R-S(\Omega_2)R,
\end{gather}
where $m=\frac{m_1m_2}{m_1+m_2}$. The following equations can be
used to reconstruct the motion of the second body in the inertial
frame:
\begin{gather}
\dot{\gamma}_2 = R_2 \deriv{U}{X},\\
\dot{x}_2 = \frac{\gamma_2}{m_2},\\
\dot{R_2}=R_2 S(\Omega_2)\label{eqn:R2doth}.
\end{gather}

\section{Variational integrators}
A variational integrator discretizes Hamilton's principle rather
than the continuous equations of motion. Taking variations of the
discretization of the action integral leads to the discrete
Euler-Lagrange or discrete Hamilton's equations. The discrete
Euler-Lagrange equations can be interpreted as a discrete Lagrangian
map that updates the variables in the configuration space, which are
the positions and the attitudes of the bodies. A discrete Legendre
transformation relates the configuration variables with the linear
and angular momenta variables, and yields a discrete Hamiltonian
map, which is equivalent to the discrete Lagrangian map.

In this section, we derive both a Lagrangian and Hamiltonian form of
variational integrators for the full body problem in inertial and
relative coordinates. The second level subscript $k$ denotes the
value of variables at $t=kh+t_0$ for an integration step size
$h\in\Rset$ and an integer $k$. The integer $N$ satisfies
$t_f=kN+t_0$, so $N$ is the number of time-steps of length $h$ to go
from the initial time $t_0$ to the final time $t_f$.

\subsection{Inertial coordinates}

\textit{Discrete Lagrangian:} In continuous time, the structure of
the kinematics equations \refeqn{Ridot}, \refeqn{Rdot} and
\refeqn{R2dot} ensure that $R_i$, $R$ and $R_2$ evolve on $\SO$
automatically. Here, we introduce a new variable $F_{i_{k}}\in\SO$
defined such that $R_{i_{k+1}}=R_{i_{k}}F_{i_{k}}$, i.e.
\begin{align}
F_{i_{k}}=R_{i_{k}}^TR_{i_{k+1}}.\label{eqn:Fik}
\end{align}
Thus, $F_{i_{k}}$ represents the relative attitude between two
integration steps, and by requiring that $F_{i_{k}}\in\SO$, we
guarantee that $R_{i_{k}}$ evolves in $\SO$.

Using the kinematic equation $\dot{R}_i=R_iS(\Omega_i)$, the
skew-symmetric matrix $S(\Omega_{_{k}})$ can be approximated as
\begin{align}
S(\Omega_{i_{k}}) & = R_{i_{k}}^T \dot{R}_{i_{k}}\approx R_{i_{k}}^T
\frac{R_{i_{k+1}}-R_{i_{k}}}{h}= \frac{1}{h}(F_{i_{k}}-I_{3\times
3}).\label{eqn:SOmegaik}
\end{align}
The velocity, $\dot{x}_{i_{k}}$ can be approximated simply by
$(x_{i_{k+1}}-x_{i_{k}})/h$. Using these approximations of the
angular and linear velocity, the kinetic energy of the $i$th body
given in \refeqn{Ti} can be approximated as
\begin{align*}
T_i(\dot{x}_i,\Omega_i)& \approx
T_i\parenth{\frac{1}{h}(x_{i_{k+1}}-x_{i_{k}}),\frac{1}{h}(F_{i_{k}}-I_{3\times 3})},\\
& =\frac{1}{2h^2}m_i\norm{x_{i_{k+1}}-x_{i_{k}}}^2
+\frac{1}{2h^2}\tr{(F_{i_{k}}-I_{3\times 3})J_{d_i}(F_{i_{k}}-I_{3\times 3})^T},\\
& =\frac{1}{2h^2}m_i\norm{x_{i_{k+1}}-x_{i_{k}}}^2
+\frac{1}{h^2}\tr{(I_{3\times 3}-F_{i_{k}})J_{d_i}},
\end{align*}
where \refeqn{trAB} is used. A discrete Lagrangian
$L_d(\mathbf{x}_{_{k}},\mathbf{x}_{_{k+1}},\mathbf{R}_{_{k}},\mathbf{F}_{_{k}})$
is constructed such that it approximates a segment of the action
integral \refeqn{caction},
\begin{align}
L_d & = \frac{h}{2}L\parenth{\mathbf{x}_{_{k}},\frac{1}{h}(\mathbf
x_{_{k+1}}-\mathbf x_{_{k}})
,\mathbf{R}_{_{k}},\frac{1}{h}(\mathbf F_{_{k}}-\mathbf
I)}\nonumber\\
& \quad +
\frac{h}{2}L\parenth{\mathbf{x}_{_{k+1}},\frac{1}{h}(\mathbf
x_{_{k+1}}-\mathbf x_{_{k}})
,\mathbf{R}_{_{k+1}},\frac{1}{h}(\mathbf F_{_{k}}-\mathbf
I)}\nonumber,\\
& = \sum_{i=1}^{n} \frac{1}{2h}m_i\norm{x_{i_{k+1}}-x_{i_{k}}}^2
+\frac{1}{h}\tr{(I_{3\times 3}-F_{i_{k}})J_{d_i}}\nonumber\\
& \quad -\frac{h}{2}U(\mathbf{x}_{_{k}},\mathbf{R}_{_{k}})
-\frac{h}{2}U(\mathbf{x}_{_{k+1}},\mathbf{R}_{_{k+1}}),\label{eqn:Ld}
\end{align}
where $\mathbf{x}_{_k}\in(\Rset^3)^n$, $\mathbf{R}_{_k}\in\SO^n$,
and $\mathbf{F}_{_k}\in(\Rset^3)^n$, and
$\mathbf{I}\in(\Rset^{3\times 3})^n$ are defined as
$\mathbf{x}_{_k}=(x_{1_k},x_{2_k},\cdots,x_{n_k})$,
$\mathbf{R}_{_k}=(R_{1_k},R_{2_k},\cdots,R_{n_k})$,
$\mathbf{F}_{_k}=(F_{1_k},F_{2_k},\cdots,F_{n_k})$, and
$\mathbf{I}=(I_{3\times 3},I_{3\times 3},\cdots,I_{3\times 3})$,
respectively.

This discrete Lagrangian is self-adjoint~\cite{bk:hairer}, and
self-adjoint numerical integration methods have even order, so we
are guaranteed that the resulting integration method is at least
second-order accurate.

\textit{Variations of discrete variables:} The variations of the
discrete variables are chosen to respect the geometry of the
configuration space $\SE$. The variation of $x_{i_{k}}$ is given by
\begin{align*}
x_{i_{k}}^\epsilon = x_{i_{k}} + \epsilon\delta x_{i_{k}}
+\mathcal{O}(\epsilon^2),
\end{align*}
where $\delta x_{i_{k}}\in\Rset^3$ and vanishes at $k=0$ and $k=N$.
The variation of $R_{i_{k}}$ is given by
\begin{align}
\delta R_{i_{k}}=R_{i_{k}}\eta_{i_{k}},\label{eqn:delRik}
\end{align}
where $\eta_{i_{k}}\in\so$ is a variation represented by a
skew-symmetric matrix and vanishes at $k=0$ and $k=N$. The variation
of $F_{i_{k}}$ can be computed from the definition
$F_{i_{k}}=R_{i_{k}}^TR_{i_{k+1}}$ to give
\begin{align}
\delta F_{i_{k}} & = \delta R_{i_{k}}^TR_{i_{k+1}} +
R_{i_{k}}^T\delta R_{i_{k+1}},\nonumber\\
& = -\eta_{i_{k}} R_{i_{k}}^TR_{i_{k+1}} + R_{i_{k}}^TR_{i_{k+1}}
\eta_{i_{k+1}},\nonumber\\
& =
-\eta_{i_{k}}F_{i_{k}}+F_{i_{k}}\eta_{i_{k+1}}.\label{eqn:delFik}
\end{align}

\subsubsection{Discrete equations of motion: Lagrangian
form}
To obtain the discrete equations of motion in Lagrangian form, we
compute the variation of the discrete Lagrangian from \refeqn{delU},
\refeqn{delRik} and \refeqn{delFik}, to give
\begin{align}
\delta L_d = \sum_{i=1}^{n}\, & \frac{1}{h} m_i
(x_{i_{k+1}}-x_{i_{k}})^T(\delta x_{i_{k+1}}-\delta x_{i_{k}})
+\frac{1}{h}\tr{\parenth{\eta_{i_{k}}F_{i_{k}}-F_{i_{k}}\eta_{i_{k+1}}}J_{d_i}}\nonumber\\
& -\frac{h}{2}\parenth{\deriv{U_{_{k}}}{x_{i_{k}}}^T\delta
x_{i_{k}}+\deriv{U_{_{k+1}}}{x_{i_{k+1}}}^T\delta x_{i_{k+1}}}
+\frac{h}{2}\tr{\eta_{i_{k}}R_{i_{k}}^T\deriv{U_{_{k}}}{R_{i_{k}}}
+\eta_{i_{k+1}}R_{i_{k+1}}^T\deriv{U_{_{k+1}}}{R_{i_{k+1}}}},\label{eqn:delLd}
\end{align}
where $U_{_{k}}=U(\mathbf{x}_{_{k}},\mathbf{R}_{_{k}})$ denotes the
value of the potential at $t=kh+t_0$.

Define the action sum as
\begin{align}
\mathfrak{G}_d =
\sum_{k=0}^{N-1}L_d(\mathbf{x}_{_{k}},\mathbf{x}_{_{k+1}},\mathbf{R}_{_{k}},\mathbf{F}_{_{k}}).\label{eqn:daction}
\end{align}
The discrete action sum $\mathfrak{G}_d$ approximates the action
integral \refeqn{caction}, because the discrete Lagrangian
approximates a segment of the action integral.

Substituting \refeqn{delLd} into \refeqn{daction}, the variation of
the action sum is given by
\begin{align*}
\delta\mathfrak{G}_d = \sum_{k=0}^{N-1}\sum_{i=1}^{n}\; &
\delta x_{i_{k+1}}^T \braces{\frac{1}{h} m_i
(x_{i_{k+1}}-x_{i_{k}})-\frac{h}{2}\deriv{U_{_{k+1}}}{x_{i_{k+1}}}}\\
& +\delta x_{i_{k}}^T \braces{-\frac{1}{h} m_i
(x_{i_{k+1}}-x_{i_{k}})-\frac{h}{2}\deriv{U_{_{k}}}{x_{i_{k}}}}\\
& +\tr{\eta_{i_{k+1}}\braces{-\frac{1}{h}J_{d_i}F_{i_k}
+\frac{h}{2}R_{i_{k+1}}^T\deriv{U_{_{k+1}}}{R_{i_{k+1}}}}}\\
& +\tr{\eta_{i_{k}}\braces{\frac{1}{h}F_{i_k}J_{d_i}
 +\frac{h}{2}R_{i_{k}}^T\deriv{U_{_{k}}}{R_{i_{k}}}}}.
\end{align*}
Using the fact that $\delta x_{i_{k}}$ and $\eta_{i_{k}}$ vanish at
$k=0$ and $k=N$, we can reindex the summation, which is the discrete
analogue of integration by parts, to yield
\begin{align*}
\delta \mathfrak{G}_d = \sum_{k=1}^{N-1}\sum_{i=1}^n\,&
-\delta x_{i_{k}}
\braces{\frac{1}{h}m_i\parenth{x_{i_{k+1}}-2x_{i_{k}}+x_{i_{k-1}}}
+h\deriv{U_{_{k}}}{x_{i_{k}}}}\\
&+\tr{\eta_{i_{k}}\braces{\frac{1}{h}\parenth{F_{i_{k}}J_{d_i}-J_{d_i}F_{i_{k-1}}}
+h R_{i_{k}}^T\deriv{U_{_{k}}}{R_{i_{k}}}}}.
\end{align*}

Hamilton's principle states that $\delta\mathfrak{G}_d$ should be
zero for all possible variations $\delta x_{i_{k}}\in\Rset^3$ and
$\eta_{i_{k}}\in\so$ that vanish at the endpoints. Therefore, the
expression in the first brace should be zero, and since
$\eta_{i_{k}}$ is skew-symmetric, the expression in the second brace
should be symmetric. Thus, we obtain \textit{the discrete equations
of motion for the full body problem, in Lagrangian form,} for
$i\in(1,2,\cdots,n)$ as
\begin{gather}
\frac{1}{h}\parenth{x_{i_{k+1}}-2x_{i_{k}}+x_{i_{k-1}}}=-h\deriv{U_{_{k}}}{x_{i_{k}}},\label{eqn:updatexik}\\
\frac{1}{h}\parenth{F_{i_{k+1}}J_{d_i}-J_{d_i}F_{i_{k+1}}^T
-J_{d_i}F_{i_{k}}+F_{i_{k}}^TJ_{d_i}}
=h S(M_{i_{k+1}}),\label{eqn:findFik}\\
R_{i_{k+1}}=R_{i_{k}}F_{i_{k}},\label{eqn:updateRik}
\end{gather}
where $M_{i_{k}}\in\Rset^3$ is defined in \refeqn{Mi} as
\begin{align}
M_{i_k}=r_{i_1}\times u_{ri_1}+r_{i_2}\times u_{ri_2}+r_{i_3}\times
u_{ri_3},\label{eqn:Mik}
\end{align}
where $r_{i_p},u_{ri_p}\in\Rset^{1\times 3}$ are $p$th row vectors
of $R_{i_k}$ and $\deriv{U_{_k}}{R_{i_k}}$, respectively. Given the
initial conditions $(x_{i_{0}},R_{i_{0}},x_{i_{1}},R_{i_{1}})$, we
can obtain $x_{i_{2}}$ from \refeqn{updatexik}. Then, $F_{i_{0}}$ is
computed from \refeqn{updateRik}, and $F_{i_{1}}$ can be obtained by
solving the implicit equation \refeqn{findFik}. Finally, $R_{i_{2}}$
is found from \refeqn{updateRik}. This yields an update map
$(x_{i_{0}},R_{i_{0}},x_{i_{1}},R_{i_{1}})\mapsto(x_{i_{1}},R_{i_{1}},x_{i_{2}},R_{i_{2}})$,
and this process can be repeated.

\vspace{-0.2cm}
\subsubsection{Discrete equations of motion: Hamiltonian
form} \vspace{-0.2cm}

As discussed above, equations \refeqn{updatexik} through
\refeqn{updateRik} defines a discrete Lagrangian map that updates
$x_{i_{k}}$ and $R_{i_{k}}$. The discrete Legendre transformation
given in \refeqn{dLtk} and \refeqn{dLtkp} relates the configuration
variables $x_{i_{k}}$, $R_{i_{k}}$ and the corresponding momenta.
This induces a discrete Hamiltonian map that is equivalent to the
discrete Lagrangian map. The discrete Hamiltonian map is
particularly convenient if the initial conditions are given in terms
of the positions and momenta at the initial time,
$(x_{i_{0}},v_{i_{0}},R_{i_{0}},\Omega_{i_{0}})$.

Before deriving the variational integrator in Hamiltonian form,
consider the momenta conjugate to $x_i$ and $R_i$, namely
$P_{v_i}\in\Rset^3$ and $P_{\Omega_i}\in\Rset^{3\times 3}$. From the
definition \refeqn{FL}, $\mathbb{F}_{v_i}L$ is obtained by taking
the derivative of $L$, given in \refeqn{L}, with $\dot{x}_i$ while
holding other variables fixed.
\begin{align*}
\delta\dot{x}_i^TP_{v_i}&=\mathbb{F}_{v_i}L(\mathbf{x},\mathbf{\dot
x},\mathbf{R},\mathbf{\Omega}),\\
&=\frac{d}{d\epsilon}\bigg|_{\epsilon=0}L(\mathbf{x},\mathbf{\dot
x}+\epsilon\delta\mathbf{\dot{x}}_i,\mathbf{R},\mathbf{\Omega}),\\
& = \frac{d}{d\epsilon}\bigg|_{\epsilon=0}
T_i(\dot{x}_i+\epsilon\delta\dot{x}_i,\Omega_i),\\
& = \delta\dot{x}_i^T \parenth{m_i \dot{x}_i},
\end{align*}
where $\delta\mathbf{\dot{x}}_i\in(\Rset^3)^n$ denotes
$(0,0,\cdots,\delta\dot{x}_i,\cdots,0)$, and $T_i$ is given in
\refeqn{Ti}. Then, we obtain
\begin{align}
P_{v_i} & = m_i v_i = \gamma_i, \label{eqn:Pvi}
\end{align}
which is equal to the linear momentum of $\mathcal{B}_i$. Similarly,
\begin{align*}
\tr{S(\delta\Omega_i)^TP_{\Omega_i}}
& =\mathbb{F}_{\Omega_i}L(\mathbf{x},\mathbf{\dot
x},\mathbf{R},\mathbf{\Omega}),\\
& = \frac{d}{d\epsilon}\bigg|_{\epsilon=0}
T_i(\dot{x}_i,\Omega_i+\epsilon\delta\Omega_i),\\
& =\frac{1}{2}\tr{S(\delta\Omega_i)J_{d_i}S(\Omega_i)^T
+S(\Omega_i)J_{d_i}S(\delta\Omega_i)^T},\\
& = \frac{1}{2}\tr{S(\delta\Omega_i)^TS(J_i\Omega_i)},
\end{align*}
where \refeqn{trAB} and \refeqn{JdJ} are used. Now, we obtain
\begin{align*}
\tr{S(\delta\Omega_i)^T\braces{P_{\Omega_i}-\frac{1}{2}S(J_i\Omega_i)}}=0.
\end{align*}
Since $S(\Omega_i)$ is skew-symmetric, the expression in the braces
should be symmetric. This implies that
\begin{align}
P_{\Omega_i}-P_{\Omega_i}^T& =S(J_i\Omega_i) =
S(\Pi_i).\label{eqn:POmegai}
\end{align}

Equations \refeqn{Pvi} and \refeqn{POmegai} give expressions for the
momenta conjugate to $x_i$ and $R_i$. Consider the discrete Legendre
transformations given in \refeqn{dLtk} and \refeqn{dLtkp}. Then,
\begin{align}
\delta
x_{i_{k}}^T\mathbf{D}_{x_{i,k}}L_d(\mathbf{x}_{_{k}},\mathbf{x}_{_{k+1}},\mathbf{R}_{_{k}},\mathbf{F}_{_{k}})&=
\frac{d}{d\epsilon}\bigg|_{\epsilon=0}L_d(\mathbf{x}_{_{k}}+\epsilon\delta\mathbf{x}_{i_k},\mathbf{x}_{_{k+1}},\mathbf{R}_{_{k}},\mathbf{F}_{_{k}}),\nonumber\\
& = -\delta x_{i_{k}}^T \bracket{\frac{1}{h} m_i
(x_{i_{k+1}}-x_{i_{k}})
+\frac{h}{2}\deriv{U_{_{k}}}{x_{i_{k}}}},\label{eqn:DxikLd0}
\end{align}
where $\delta\mathbf{x_{i_k}}\in(\Rset^3)^n$ denotes
$(0,0,\cdots,\delta x_{_{i_k}},\cdots,0)$. Therefore, we have
\begin{align}
\mathbf{D}_{x_{i,k}}L_d(\mathbf{x}_{_{k}},\mathbf{x}_{_{k+1}},\mathbf{R}_{_{k}},\mathbf{F}_{_{k}})&=
-\frac{1}{h} m_i (x_{i_{k+1}}-x_{i_{k}})
-\frac{h}{2}\deriv{U_{_{k}}}{x_{i_{k}}}.\label{eqn:DxikLd}
\end{align}

From the discrete Legendre transformation given in \refeqn{dLtk},
$P_{v_{i,k}}=-\mathbf{D}_{x_{i,k}}L_d$. Using \refeqn{Pvi} and
\refeqn{DxikLd}, we obtain
\begin{align}
\gamma_{i_k}=\frac{1}{h} m_i (x_{i_{k+1}}-x_{i_{k}})
+\frac{h}{2}\deriv{U_{_{k}}}{x_{i_{k}}}.\label{eqn:mvk}
\end{align}
Using the discrete Legendre transformation given in \refeqn{dLtkp},
$P_{v_{i,k+1}}=\mathbf{D}_{x_{i,k+1}}L_d$, we can derive the
following equation similarly:
\begin{align}
\gamma_{i_{k+1}}=\frac{1}{h} m_i (x_{i_{k+1}}-x_{i_{k}})
-\frac{h}{2}\deriv{U_{_{k+1}}}{x_{i_{k+1}}}.\label{eqn:mvkp}
\end{align}

Equations \refeqn{mvk} and \refeqn{mvkp} define the variational
integrator in Hamiltonian form for the translational motion. Now,
consider the rotational motion. We have
\begin{align}
\tr{\eta_{i_{k}}{\mathbf{D}_{R_{i,k}}L_d}^T}&=
\tr{\eta_{i_{k}}\braces{\frac{1}{h}F_{i_{k}}J_{d_i}
+\frac{h}{2}R_{i_{k}}^T\deriv{U_{_{k}}}{R_{i_{k}}}}},\label{eqn:DRikLd0}
\end{align}
where the right side is obtained by taking the variation of $L_d$
with respect to $R_{i_{k}}$, while holding other variables fixed.
Since $\eta_{i_{k}}$ is skew-symmetric,
\begin{align}
-\mathbf{D}_{R_{i,k}}L_d+{\mathbf{D}_{R_{i,k}}L_d}^T =
\frac{1}{h}\parenth{F_{i_{k}}J_{d_i}-J_{d_i}F_{i_{k}}^T}
-\frac{h}{2}S(M_{i_k}),\label{eqn:DRikLd}
\end{align}
where $M_{i_i}\in\Rset^3$ is defined in \refeqn{Mik}.

From the discrete Legendre transformation given in \refeqn{dLtk},
$P_{\Omega_{i,k}}=-\mathbf{D}_{R_{i,k}}L_d$, we obtain the following
equation by using \refeqn{POmegai} and \refeqn{DRikLd},
\begin{align}
S(\Pi_{i_k})=\frac{1}{h}\parenth{F_{i_{k}}J_{d_i}-J_{d_i}F_{i_{k}}^T}
-\frac{h}{2}S(M_{i_k}).\label{eqn:SJOmegaik}
\end{align}
Using the discrete Legendre transformation given in \refeqn{dLtkp},
$P_{\Omega_{i,k+1}}=\mathbf{D}_{R_{i,k+1}}L_d$, we can obtain the
following equation:
\begin{align}
S(\Pi_{i_{k+1}})=\frac{1}{h}F_{i_{k}}^T\parenth{F_{i_{k}}J_{d_i}-J_{d_i}F_{i_{k}}^T}F_{i_{k}}
+\frac{h}{2}S(M_{i_{k+1}}).\label{eqn:SJOmegaikp}
\end{align}
By using \refeqn{SR} and substituting \refeqn{SJOmegaik}, we can
reduce \refeqn{SJOmegaikp} to the following equation in vector form.
\begin{align}
\Pi_{i_{k+1}}=F_{i_{k}}^T\Pi_{i_k}+\frac{h}{2}F_{i_{k}}^TM_{i_k}+\frac{h}{2}M_{i_{k+1}}.
\label{eqn:JOmegaikp}
\end{align}
Equations \refeqn{SJOmegaik} and \refeqn{JOmegaikp} define the
variational integrator in Hamiltonian form for the rotational
motion.

In summary, using \refeqn{mvk}, \refeqn{mvkp}, \refeqn{SJOmegaik}
and \refeqn{JOmegaikp}, \textit{the discrete equations of motion for
the full body problem, in Hamiltonian form,} can be written for
$i\in(1,2,\cdots,n)$ as
\begin{gather}
x_{i_{k+1}} = x_{i_{k}} + \frac{h}{m_i}
\gamma_{i_k}-\frac{h^2}{2m_i}\deriv{U_{_{k}}}{x_{i_{k}}},\label{eqn:updatexikH}\\
 \gamma_{i_{k+1}}=\gamma_{i_k}-\frac{h}{2}\deriv{U_{_{k}}}{x_{i_{k}}}
-\frac{h}{2}\deriv{U_{_{k+1}}}{x_{i_{k+1}}},\label{eqn:updatevik}\\
h
S(\Pi_{i_k}+\frac{h}{2}M_{i_k})=F_{i_{k}}J_{d_i}-J_{d_i}F_{i_{k}}^T,\label{eqn:findFikH}\\
\Pi_{i_{k+1}}=F_{i_{k}}^T\Pi_{i_k}+\frac{h}{2}F_{i_{k}}^TM_{i_k}+\frac{h}{2}M_{i_{k+1}}.
\label{eqn:updatePiik},\\
R_{i_{k+1}}=R_{i_{k}}F_{i_{k}}.\label{eqn:updateRikH}
\end{gather}
Given $(x_{i_0},\gamma_{i_0},R_{i_0},\Pi_{i_0})$, we can find
$x_{i_1}$from \refeqn{updatexikH}. Solving the implicit equation
\refeqn{findFikH} yields $F_{i_0}$, and $R_{i_1}$ is computed from
\refeqn{updateRikH}. Then,  \refeqn{updatevik} and
\refeqn{updatePiik} gives $\gamma_{i_1}$, and $\Pi_{i_1}$. This
defines the discrete Hamiltonian map,
$(x_{i_0},\gamma_{i_0},R_{i_0},\Pi_{i_0})\mapsto(x_{i_1},\gamma_{i_1},R_{i_1},\Pi_{i_1})$,
and this process can be repeated.

\vspace{-0.3cm}
\subsection{Relative coordinates}
\vspace{-0.5cm} In this section, we derive the variational
integrator for the full two body problem in relative coordinates by
following the procedure given before. This result can be readily
generalized to $n$ bodies.

\textit{Reduction of discrete variables:} The discrete reduced
variables are defined in the same way as the continuous reduced
variables, which are given in \refeqn{X} through \refeqn{V2}. We
introduce $F_{_k}\in\SO$ such that
$R_{_{k+1}}=R_{2_{k+1}}^TR_{1_{k+1}}=F_{2_k}^T F_{_k} R_{_k}$. i.e.
\begin{align}
F_{_k} = R_{_k} F_{1_k} R_{_k}^T.\label{eqn:Fk}
\end{align}

\textit{Discrete reduced Lagrangian:} The discrete reduced
Lagrangian is obtained by expressing the original discrete
Lagrangian given in \refeqn{Ld} in terms of the discrete reduced
variables.

From the definition of the discrete reduced variables given in
\refeqn{X} and \refeqn{X2}, we have
\begin{align}
x_{1_{k+1}}-x_{1_{k}}&=R_{2_{k+1}}(X_{_{k+1}}+X_{2_{k+1}})-R_{2_{k}}(X_{_{k}}+X_{2_{k}}),\nonumber\\
&=R_{2_{k}}\braces{F_{2_{k}}(X_{_{k+1}}+X_{2_{k+1}})-(X_{_{k}}+X_{2_{k}})},\label{eqn:delx1k}\\
x_{2_{k+1}}-x_{2_{k}}&=R_{2_{k}}\braces{F_{2_{k}}X_{2_{k+1}}-X_{2_k}}.\label{eqn:delx2k}
\end{align}
From \refeqn{SOmegaik}, $S(\Omega_{1_k})$ and $S(\Omega_{2_k})$ are
expressed as
\begin{align}
S(\Omega_{1_k})&=\frac{1}{h}\parenth{F_{1_{k}}-I_{3\times 3}},\nonumber\\
& =\frac{1}{h}R_{_k}^T \parenth{F_{_k} - I_{3\times 3}}
R_{_k},\label{eqn:SOmega1k}\\
S(\Omega_{2_{k}}) & = \frac{1}{h}\parenth{F_{2_{k}}-I_{3\times
3}}.\label{eqn:SOmega2k}
\end{align}
Substituting \refeqn{delx1k} through \refeqn{SOmega2k} into
\refeqn{Ld}, we obtain the discrete reduced Lagrangian.
\begin{align}
l_{d_k}& =
l_d(X_{_{k}},X_{_{k+1}},X_{2_{k}},X_{2_{k+1}},R_{_{k}},F_{_{k}},F_{2_{k}})\nonumber\\
& =
\frac{1}{2h}m_1\norm{F_{2_{k}}(X_{_{k+1}}+X_{2_{k+1}})-(X_{_k}+X_{2_{k}})}^2
+\frac{1}{2h}m_2\norm{F_{2_{k}}X_{2_{k+1}}-X_{2_k}}^2\nonumber\\
& \quad + \frac{1}{h}\tr{(I_{3\times 3}-F_{_{k}})J_{{dR}_k}} +
\frac{1}{h}\tr{(I_{3\times
3}-F_{2_{k}})J_{d_2}}-\frac{h}{2}U(X_{_{k}},R_{_{k}})-\frac{h}{2}U(X_{_{k+1}},R_{_{k+1}}),\label{eqn:ld}
\end{align}
where $J_{{dR}_k}\in\Rset^{3\times 3}$ is defined to be
$J_{{dR}_k}=R_{_{k}}J_{d_1}R_{_{k}}^T$, which gives the nonstandard
moment of inertia matrix of the first body with respect to the
second body fixed frame at $t=kh+t_0$.

\textit{Variations of discrete reduced variables:} The variations of
the discrete reduced variables can be derived from those of the
original variables. The variations of $R_{_{k}}$, $X_{_{k}}$, and
$F_{2_{k}}$ are the same as given in \refeqn{delR}, \refeqn{delX},
and \refeqn{delFik}, respectively. The variation of $F_{_{k}}$ is
computed in \ref{appdrv}.

In summary, the variations of discrete reduced variables are given
by
\begin{align}
\delta R_{_k} & = \eta_{_k} R_{_k} - \eta_{2_k}
R_{_k},\label{eqn:delRk}\\
\delta X_{_k} & = \chi_{_k} -\eta_{2_k}X_{_k},\label{eqn:delXk}\\
\delta F_{_k} & =-\eta_{2_{k}}F_{_{k}}
+F_{2_{k}}\eta_{_{k+1}}F_{2_{k}}^TF_{_{k}}
+F_{_{k}}\parenth{-\eta_{_{k}} +\eta_{2_{k}} },\label{eqn:delFk}\\
\delta X_{2_k} & = \chi_{2_k} -\eta_{2_k}X_{2_k},\label{eqn:delX2k}\\
\delta F_{2_k} &
=-\eta_{2_k}F_{2_k}+F_{2_k}\eta_{2_{k+1}}.\label{eqn:delF2k}
\end{align}
These Lie group variations are the main elements required to derive
the variational integrator equations.

\subsubsection{Discrete equations of motion: Lagrangian
form} As before, we can obtain the discrete equations of motion in
Lagrangian form by computing the variation of the discrete reduced
Lagrangian which, by using \refeqn{delRk} through \refeqn{delF2k},
is given by
\begin{align}
\delta l_{d_k} & =
\frac{1}{h}\chi_{_{k+1}}^T\bracket{m_1(X_{_{k+1}}+X_{2_{k+1}})-m_1F_{2_{k}}^T(X_{_k}+X_{2_{k}})}\nonumber\\
& \quad
+\frac{1}{h}\chi_{_{k}}^T\bracket{m_1(X_{_{k}}+X_{2_{k}})-m_1F_{2_{k}}(X_{_{k+1}}
+X_{2_{k+1}})}\nonumber\\
& \quad +
\frac{1}{h}\chi_{2_{k+1}}^T\bracket{m_1(X_{_{k+1}}+X_{2_{k+1}})-m_1F_{2_{k}}^T(X_{_k}+X_{2_{k}})
+m_2X_{2_{k+1}}-m_2F_{2_{k}}^TX_{2_{k}}}\nonumber\\
& \quad + \frac{1}{h}
\chi_{2_{k}}^T\bracket{m_1(X_{_{k}}+X_{2_{k}})-m_1F_{2_{k}}(X_{_{k+1}}+X_{2_{k+1}})
+m_2X_{2_{k}}-m_2F_{2_{k}}X_{2_{k+1}}}\nonumber\\
& \quad -\frac{1}{h}\tr{\eta_{_{k+1}}F_{2_{k}}^TF_{_{k}} J_{{dR}_k}
F_{2_{k}}}
+\frac{1}{h}\tr{\eta_{_{k}}F_{_{k}}J_{{dR}_k}}
-\frac{1}{h}\tr{\eta_{2_{k+1}}
J_{d_2}F_{2_{k}}}+\frac{1}{h}\tr{\eta_{2_{k}}F_{2_{k}}J_{d_2}}\nonumber\\
& \quad - \frac{h}{2}\chi_{_{k}}^T\deriv{U_{_k}}{X_{_{k}}}
+\frac{h}{2}\tr{\eta_{2_{k}}X_{_{k}}\deriv{U_{_k}}{X_{_{k}}}^T}
- \frac{h}{2}\chi_{_{k+1}}^T\deriv{U_{_{k+1}}}{X_{_{k+1}}}
+\frac{h}{2}\tr{\eta_{2_{k+1}}X_{_{k+1}}\deriv{U_{_{k+1}}}{X_{_{k+1}}}^T}\nonumber\\
& \quad
+\frac{h}{2}\tr{\eta_{2_{k}}R_{_{k}}\deriv{U_{_k}}{R_{_{k}}}^T
-\eta_{_{k}}R_{_{k}}\deriv{U_{_k}}{R_{_{k}}}^T}
+\frac{h}{2}\tr{\eta_{2_{k+1}}R_{_{k+1}}\deriv{U_{_{k+1}}}{R_{_{k+1}}}^T
-\eta_{_{k+1}}R_{_{k+1}}\deriv{U_{_{k+1}}}{R_{_{k+1}}}^T}\label{eqn:delldk}.
\end{align}

The action sum expressed in terms of the discrete reduced Lagrangian
has the form
\begin{align}
\mathfrak{G}_d =
\sum_{k=0}^{N-1}l_d(X_{_{k}},X_{_{k+1}},X_{2_{k}},X_{2_{k+1}},R_{_{k}},F_{_{k}},F_{2_{k}}).\label{eqn:draction}
\end{align}
The discrete action sum $\mathfrak{G}_d$ approximates the action
integral \refeqn{craction}, because the discrete Lagrangian
approximates a piece of the integral. Using the fact that the
variations $\chi_{_{k}},\chi_{2_{k}},\eta_{_{k}},\eta_{2_{k}}$
vanish at $k=0$ and $k=N$, the variation of the discrete action sum
can be expressed as
\begin{align*}
\delta \mathfrak{G}_d = & \sum_{k=1}^{N-1}
\frac{1}{h}\chi_{_{k}}^T\bigg\{-m_1F_{2_{k-1}}^T(X_{_{k-1}}+X_{2_{k-1}})
+2m_1(X_{_{k}}+X_{2_{k}})\\
& \hspace{1.8cm} -m_1F_{2_{k}}(X_{_{k+1}}
+X_{2_{k+1}})-h^2\deriv{U_{_k}}{X_{_{k}}}\bigg\}\\
& + \sum_{k=1}^{N-1}
\frac{1}{h}\chi_{2_{k}}^T\bigg\{-m_1F_{2_{k-1}}^T(X_{_{k-1}}+X_{2_{k-1}})
+2m_1(X_{_{k}}+X_{2_{k}})-m_1F_{2_{k}}(X_{_{k+1}} +X_{2_{k+1}})\\
& \hspace{2.2cm}
-m_2F_{2_{k-1}}^TX_{2_{k-1}}+2m_2X_{2_{k}}-m_2F_{2_{k}}X_{2_{k+1}}\bigg\}\\
& + \sum_{k=1}^{N-1}
\tr{\eta_{_{k}}\braces{\frac{1}{h}\parenth{-F_{2_{k-1}}^TF_{_{k-1}}
R_{_{k-1}}J_{d_1}R_{_{k-1}}^T F_{2_{k-1}}
+F_{_{k}}R_{_{k}}J_{d_1}R_{_{k}}^T}-hR_{_{k}}\deriv{U_{_k}}{R_{_{k}}}^T}}\\
& + \sum_{k=1}^{N-1} \tr{\eta_{2_{k}}
\braces{\frac{1}{h}\parenth{-J_{d_2}F_{2_{k-1}}+F_{2_{k}}J_{d_2}}
+hX_{_{k}}\deriv{U_{_k}}{X_{_{k}}}^T+hR_{_{k}}\deriv{U_{_k}}{R_{_{k}}}^T}}.
\end{align*}
From Hamilton's principle, $\delta \mathfrak{G}_d$ should be zero for
all possible variations $\chi_{_{k}},\chi_{2_{k}}\in\Rset^3$ and
$\eta_{_{k}},\eta_{2_{k}}\in\so$ which vanish at the endpoints.
Therefore, in the above equation, the expressions in the first two
braces should be zero, and the expressions in the last two braces
should be symmetric since $\eta_{_{k}},\eta_{2_{k}}$ are
skew-symmetric. After some simplification, we obtain \textit{the
discrete equations of relative motion for the full two body problem,
in Lagrangian form,} as
\begin{gather}
F_{2_{k}}X_{_{k+1}}-2X_{_{k}}+F_{2_{k-1}}^TX_{_{k-1}}=-\frac{h^2}{m}\deriv{U_{_k}}{X_{_{k}}},\label{eqn:updateXl}\\
F_{_{k+1}}J_{dR_{k+1}}-J_{dR_{k+1}}F_{_{k+1}}^T
=F_{2_{k}}^T
\parenth{F_{_{k}}J_{dR_{k}}-J_{dR_{k}}F_{_{k}}^T}F_{2_{k}}
-h^2S(M_{k+1}),\label{eqn:findFl}\\
F_{2_{k+1}}J_{d_2}-J_{d_2}F_{2_{k+1}}^T
=F_{2_{k}}^T\parenth{F_{2_{k}}J_{d_2}-J_{d_2}F_{2_{k}}^T}F_{2_{k}}
+h^2X_{_{k+1}}\times\deriv{U}{X_{_{k+1}}}+h^2S(M_{k+1})
,\label{eqn:findF2l}\\
R_{_{k+1}}=F_{2_{k}}^TF_{_{k}}R_{_{k}},\label{eqn:updateRl}\\
R_{2_{k+1}}=R_{2_{k}}F_{2_{k}}.\label{eqn:updateR2l}
\end{gather}
It is natural to express equations of motion for the second body in
the inertial frame.
\begin{gather}
x_{2_{k+1}}-2x_{2_{k}}+x_{2_{k-1}}=\frac{h^2}{m_2}R_{_k}\deriv{U_{_k}}{X_{_{k}}}.\label{eqn:updatex2l}
\end{gather}
Given $(X_{_{0}},R_{_0},R_{2_{0}},X_{_{1}},R_{_1},R_{2_{1}})$, we
can determine $F_{_0}$ and $F_{2_0}$ from \refeqn{updateRl} and
\refeqn{updateR2l}. Solving the implicit equations \refeqn{findFl}
and \refeqn{findF2l} gives $F_{_1}$ and $F_{2_1}$. Then $X_{_2}$,
$R_{_2}$ and $R_{2_2}$ are found from \refeqn{updateXl},
\refeqn{updateRl} and \refeqn{updateR2l}, respectively. This yields
the discrete Lagrangian map
$(X_{_{0}},R_{_0},R_{2_{0}},X_{_{1}},R_{_1},R_{2_{1}})\mapsto
(X_{_{1}},R_{_1},R_{2_{1}},X_{_{2}},R_{_2},R_{2_{2}})$ and this
process can be repeated. We can separately reconstruct $x_{2_k}$
using \refeqn{updatex2l}.

\vspace{-0.3cm}
\subsubsection{Discrete equations of motion: Hamiltonian
form} \vspace{-0.3cm} Using the discrete Legendre transformation, we
can obtain the Hamiltonian map, in terms of reduced variables, that
is equivalent to the Lagrangian map given in \refeqn{updateXl}
through \refeqn{updatex2l}. We will only sketch the procedure as it
is analogous to the approach of the previous section. First, we find
expressions for the conjugate momenta variables corresponding to
\refeqn{Pvi} and \refeqn{POmegai}. We compute the discrete Legendre
transformation by taking the variation of the discrete reduced
Lagrangian as in \refeqn{DxikLd0} and \refeqn{DRikLd0}. Then, we
obtain the discrete equations of motion in Hamiltonian form using
\refeqn{dLtk} and \refeqn{dLtkp}.

\textit{The discrete equations of relative motion for the full two
body problem, in Hamiltonian form,} can be written as
\begin{gather}
X_{_{k+1}}=F_{2_{k}}^T\parenth{X_{_{k}}+h\frac{\Gamma_{_{k}}}{m}-\frac{h^2}{2m}\deriv{U_{_k}}{X_{_{k}}}},
\label{eqn:updateX}\\
\Gamma_{_{k+1}}=F_{2_{k}}^T\parenth{\Gamma_{_{k}}-\frac{h}{2}\deriv{U_{_k}}{X_{_{k}}}}
-\frac{h}{2}\deriv{U_{_{k+1}}}{X_{_{k+1}}},
\label{eqn:updateV}\\
\Pi_{_{k+1}}=F_{2_{k}}^T
\parenth{\Pi_{_{k}}-\frac{h}{2}M_{_{k}}}-\frac{h}{2}M_{_{k+1}},
\label{eqn:updatePi}\\
\Pi_{2_{k+1}}=F_{2_{k}}^T\parenth{\Pi_{_{k}}+\frac{h}{2}X_{_{k}}\times\deriv{U}{X_{_{k}}}
+\frac{h}{2}M_{_{k}}}+\frac{h}{2}X_{_{k+1}}\times\deriv{U}{X_{_{k+1}}}+\frac{h}{2}M_{_{k+1}},
\label{eqn:updatePi2}\\
R_{_{k+1}}=F_{2_{k}}^TF_{_{k}}R_{_{k}},\label{eqn:updateR}\\
h S\parenth{\Pi_k-\frac{h}{2}M_{_{k}}} =
F_{_{k}}J_{dR_k}-J_{dR_k}F_{_{k}}^T,\label{eqn:findF}\\
h
S\parenth{\Pi_{2_{k}}+\frac{h}{2}X_{_{k}}\times\deriv{U}{X_{_{k}}}+\frac{h}{2}M_{_{k}}}
=F_{2_{k}}J_{d_2}-J_{d_2}F_{2_{k}}^T.\label{eqn:findF2}
\end{gather}
It is natural to express equations of motion for the second body in
the inertial frame for reconstruction:
\begin{gather}
x_{2_{k+1}}=x_{2_{k}}+h\frac{\gamma_{2_{k}}}{m_2}+\frac{h^2}{2m_2}R_{_k}\deriv{U_{_k}}{X_{_{k}}},
\label{eqn:updateX2}\\
\gamma_{2_{k+1}}=\gamma_{2_{k}}+\frac{h}{2}R_{_k}\deriv{U_{_k}}{X_{_{k}}}
+\frac{h}{2}R_{_{k+1}}\deriv{U_{_{k+1}}}{X_{_{k+1}}},
\label{eqn:updateV2}\\
R_{2_{k+1}}=R_{2_{k}}F_{2_{k}}.\label{eqn:updateR2}
\end{gather}
Given $(R_{_0},X_{_0},\Pi_{_0},\Gamma_{_0},\Pi_{2_0})$, we can
determine $F_{_0}$ and $F_{2_0}$ by solving the implicit equations
\refeqn{findF} and \refeqn{findF2}. Then, $X_{_1}$ and $R_{_1}$ are
found from \refeqn{updateX} and \refeqn{updateR}, respectively.
After that, we can compute $\Gamma_{_1}$, $\Pi_{_1}$, and
$\Pi_{2_1}$ from \refeqn{updateV}, \refeqn{updatePi} and
\refeqn{updatePi2}. This yields a discrete Hamiltonian map
$(R_{_0},X_{_0},\Pi_{_0},\Gamma_{_0},\Pi_{2_0})\mapsto
(R_{_1},X_{_1},\Pi_{_1},\Gamma_{_1},\Pi_{2_1})$, and this process
can be repeated. $x_{2_k}$, $\gamma_{2_k}$ and $R_{2_k}$ can be
updated separately using \refeqn{updateX2}, \refeqn{updateV2} and
\refeqn{updateR2}, respectively, for reconstruction.

\subsection{Numerical considerations}
\textit{Properties of the variational integrators:} Variational
integrators exhibit a discrete analogue of Noether's
theorem~\cite{jo:marsden}, and symmetries of the discrete Lagrangian
result in conservation of the corresponding momentum maps. Our
choice of discrete Lagrangian is such that it inherits the
symmetries of the continuous Lagrangian. Therefore, all the
conserved momenta in the continuous dynamics are preserved by the
discrete dynamics.

The proposed variational integrators are expressed in terms of Lie
group computations~\cite{ic:iserles}. During each integration step,
$F_{i_k} \in \SO$ is obtained by solving an implicit equation, and
$R_{i_k}$ is updated by multiplication with $F_{i_k}$. Since $\SO$
is closed under matrix multiplication, the attitude matrix
$R_{i_{k+1}}$ remains in \SO. We make this more explicit in section
\ref{comp} by expressing $F_{i_k}$ as the exponential function of an
element of the Lie algebra $\mathfrak{so}(3)$.

An adjoint integration method is the inverse map of the original
method with reversed time-step. An integration method is called
self-adjoint or symmetric if it is identical with its adjoint; a
self-adjoint method always has even order. Our discrete Lagrangian
is chosen to be self-adjoint, and therefore the corresponding
variational integrators are second-order accurate.

\textit{Higher-order methods:} While the numerical methods we
present in this paper are second-order, it is possible to apply the
symmetric composition methods, introduced in~\cite{jo:Yoshida}, to
construct higher-order versions of the Lie group variational
integrators introduced here. Given a basic numerical method
represented by the flow map $\Phi_h$, the composition method is
obtained by applying the basic method using different step sizes,
\[ \Psi_h = \Phi_{\lambda_s h}\circ\ldots\circ\Phi_{\lambda_1 h},\]
where $\lambda_1,\lambda_2,\cdots,\lambda_s\in\Rset$. In particular,
the Yoshida symmetric composition method for composing a symmetric
method of order 2 into a symmetric method of order 4 is obtained
when $s=3$, and
\[ \lambda_1=\lambda_3=\frac{1}{2-2^{1/3}},\qquad \lambda_2=-\frac{2^{1/3}}{2-2^{1/3}}.\]

Alternatively, by adopting the formalism of higher-order Lie group
variational integrators introduced in~\cite{jo:Leok2004} in
conjunction with the Rodrigues formula, one can directly construct
higher-order generalizations of the Lie group methods presented
here.

\textit{Reduction of orthogonality loss due to roundoff error:} In
the Lie group variational integrators, the numerical solution is
made to automatically remain on the rotation group by requiring that
the numerical solution is updated by matrix multiplication with the
exponential of a skew symmetric matrix.

Since the exponential of the skew symmetric matrix is orthogonal to
machine precision, the numerical solution will only deviate from
orthogonality due to the accumulation of roundoff error in the
matrix multiplication, and this orthogonality loss grows linearly
with the number of timesteps taken.

One possible method of addressing this issue is to use the
Baker-Campbell-Hausdorff (BCH) formula to track the updates purely
at the level of skew symmetric matrices (the Lie algebra). This
allows us to find a matrix $C(t)$, such that,
\[ \exp(tA)\exp(tB)=\exp C(t).\]
This matrix $C(t)$ satisfies the following differential equation,
\[ \dot C=A+B +\frac{1}{2}[A-B,C]
+\sum_{k \geq 2}\frac{B_k}{k!}\operatorname{ad}_C^k(A+B),\] with
initial value $C(0)=0$, and where $B_k$ denotes the Bernoulli
numbers, and $\operatorname{ad}_C A = [C,A]=CA-AC$.

The problem with this approach is that the matrix $C(t)$ is not
readily computable for arbitrary $A$ and $B$, and in practice, the
series is truncated, and the differential equation is solved
numerically.

An error is introduced in truncating the series, and numerical
errors are introduced in numerically integrating the differential
equations. Consequently, while the BCH formula could be used solely
at the reconstruction stage to ensure that the numerical attitude
always remains in the rotation group to machine precision, the
truncation error would destroy the symplecticity and momentum
preserving properties of the numerical scheme.

However, by combining the BCH formula with the Rodrigues formula in
constructing the discrete variational principle, it might be
possible to construct a Lie group variational integrator that tracks
the reconstructed trajectory on the rotation group at the level of a
curve in the Lie algebra, while retaining its structure-preservation
properties.

\subsection{Computational approach}\label{comp}
The structure of the discrete equations of motion given in
\refeqn{findFik}, \refeqn{findFikH}, \refeqn{findFl},
\refeqn{findF2l}, \refeqn{findF}, and \refeqn{findF2} suggests a
specific computational approach. For a given $g\in\Rset^3$, we have
to solve the following Lyapunov-like equation to find $F\in\SO$ at
each integration step.
\begin{align}
FJ_d - J_d F^T = S(g).\label{eqn:findf}
\end{align}

We now introduce an iterative approach to solve \refeqn{findf}
numerically. An element of a Lie group can be expressed as the
exponential of an element of its Lie algebra, so $F\in\SO$ can be
expressed as an exponential of $S(f)\in\so$ for some vector
$f\in\Rset^3$. The exponential can be written in closed form, using
Rodrigues' formula,
\begin{align}
F &= e^{S(f)},\nonumber\\
& = I_{3\times3} + \frac{\sin\norm{f}}{\norm{f}} S(f) +
\frac{1-\cos\norm{f}}{\norm{f}^2}S(f)^2.\label{eqn:rodc}
\end{align}
Substituting \refeqn{rodc} into \refeqn{findf}, we obtain
\begin{align*}
S(g) & = \frac{\sin\norm{f}}{\norm{f}} S(Jf) +
\frac{1-\cos\norm{f}}{\norm{f}^2} S(f\times Jf),
\end{align*}
where \refeqn{Scross} and \refeqn{JdJ} are used. Thus,
\refeqn{findf} is converted into the equivalent vector equation
$g=G(f)$, where $G : \Rset^3 \mapsto \Rset^3$ is
\begin{align*}
G(f) & = \frac{\sin\norm{f}}{\norm{f}}\, J f +
\frac{1-\cos\norm{f}}{\norm{f}^2}\, f \times J f.
\end{align*}
We use the Newton method to solve $g=G(f)$, which gives the
iteration
\begin{align}
f_{i+1} = f_i + \nabla G(f_i)^{-1} \parenth{g-
G(f_i)}.\label{eqn:newton}
\end{align}
We iterate until $\norm{g- G(f_i) } < \epsilon$ for a small
tolerance $\epsilon > 0$. The Jacobian $\nabla G(f)$ in
\refeqn{newton} can be expressed as
\begin{align*}
\nabla G(f) & =
\frac{\cos\norm{f}\norm{f}-\sin\norm{f}}{\norm{f}^3}Jff^T
+ \frac{\sin\norm{f}}{\norm{f}} J\\
& \quad +
\frac{\sin\norm{f}\norm{f}-2(1-\cos\norm{f})}{\norm{f}^4}\parenth{f
\times Jf}f^T\\
& \quad +\frac{1-\cos\norm{f}}{\norm{f}^2} \braces{-S(Jf)+S(f)J}.
\end{align*}
Numerical simulations show that 3 or 4 iterations are sufficient to
achieve a tolerance of $\epsilon=10^{-15}$.

\section{Numerical simulations}
The variational integrator in Hamiltonian form given in
\refeqn{updateX} through \refeqn{updateR2} is used to simulate the
dynamics of two simple dumbbell bodies acting under their mutual
gravity.

\subsection{Full body problem defined by two dumbbell bodies}

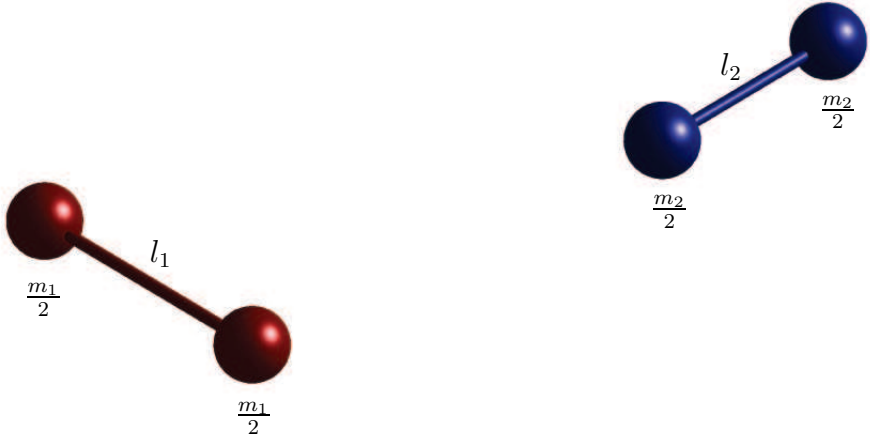
\begin{figure}[htbp]
$$\begin{xy}
\xyWARMprocessEPS{dumbbells}{eps}
\xyMarkedImport{}
\xyMarkedMathPoints{1-6}
\end{xy}$$
\caption{Dumbbell model of the full body
problem}\label{fig:dumbbell}
\end{figure}

Each dumbbell model consists of two equal rigid spheres and a
massless rod as shown in \reffig{dumbbell}.
The gravitational potential of the two dumbbell models is given by
\begin{align*}
U(X,R) =
-\sum_{p,q=1}^{2}\frac{Gm_1m_2/4}{\norm{X+\rho_{2_p}+R\rho_{1_q}}},
\end{align*}
where $G$ is the universal gravitational constant, $m_i\in\Rset$ is
the total mass of the $i$th dumbbell, and $\rho_{i_p}\in\Rset^3$ is
a vector from the origin of the body fixed frame to the $p$th sphere
of the $i$th dumbbell in the $i$th body fixed frame. The vectors
$\rho_{i_1}=[l_i/2,0,0]^T$, $\rho_{i_2}=-\rho_{i_1}$, where $l_i$ is
the length between the two spheres.

\renewcommand{\xyWARMinclude}[1]{\includegraphics[width=0.8\textwidth]{#1}}

\textit{Normalization:} Mass, length and time dimensions are
normalized as follows,
\begin{align*}
\bar{m}_i & = \frac{m_i}{m},\\
\bar{X}_i & = \frac{X_i}{l},\\
\bar{t} & = \sqrt{\frac{G(m_1+m_2)}{l^3}}\,t,
\end{align*}
where $m=\frac{m_1m_2}{m_1+m_2}$, and $l$ is chosen as the initial
horizontal distance between the center of mass of the two dumbbells.
The time is normalized so that the orbital period is of order unity.
Over-bars denote normalized variables. We can expresses the
equations of motion in terms of the normalized variables. For
example, \refeqn{Vdot} can be written as
\begin{gather*}
\bar{V}'+\bar{\Omega}_2\times \bar{V} = -\deriv{\bar{U}}{\bar{X}},
\end{gather*}
where $'$ denotes a derivative with respect to $\bar{t}$. The
normalized gravitational potential and its partial derivatives are
given by
\begin{align*}
\bar{U}&=-\frac{1}{4}\sum_{p,q=1}^{2}\frac{1}{\norm{\bar X+\bar
\rho_{2_p}+R\bar\rho_{1_q}}},\\
\deriv{\bar{U}}{\bar X}&=\frac{1}{4}\sum_{p,q=1}^{2}\frac{\bar
X+\bar \rho_{2_p}+R\bar\rho_{1_q}}{\norm{\bar X+\bar
\rho_{2_p}+R\bar\rho_{1_q}}^3},\\
\deriv{\bar{U}}{R}&=\frac{1}{4}\sum_{p,q=1}^{2}\frac{(\bar X+\bar
\rho_{2_p})\bar\rho_{1_q}^T}{\norm{\bar X+\bar
\rho_{2_p}+R\bar\rho_{1_q}}^3}.\\
\end{align*}

\textit{Conserved quantities:} The total energy $E$ is conserved:
\begin{align*}
E & =\frac{1}{2}m_1\norm{V+V_2}^2 +\frac{1}{2}m_2\norm{V_2}^2\\
& \quad +\frac{1}{2}\tr{S(\Omega)J_{d_R}S(\Omega)^T}
+\frac{1}{2}\tr{S(\Omega_2)J_{d_2}S(\Omega_2)^T}+U(X,R).
\end{align*}
 The total linear momentum $\gamma_T\in\Rset^3$, and the total angular momentum about the mass
center of the system $\pi_T\in\Rset^3$, in the inertial frame, are
also conserved:
\begin{align*}
\gamma_T & = R_2\braces{m_1\parenth{V+V_2}+m_2 V_2},\\
\pi_T & = R_2\braces{m X\times V +  J_R \Omega + J_2 \Omega_2}.
\end{align*}

\subsection{Simulation results}

The properties of the two dumbbell bodies are chosen to be
\begin{alignat*}{3}
\bar{m}_1&=1.5,&\quad \bar l_1&=0.25,&\quad \bar
J_1&=\mathrm{diag}\begin{bmatrix}0.0004&0.0238&0.0238\end{bmatrix},\\
\bar{m}_2&=3,&\quad \bar l_2&=0.5,&\quad \bar
J_2&=\mathrm{diag}\begin{bmatrix}0.0030&0.1905&0.1905\end{bmatrix}.
\end{alignat*}
The mass and length of the second dumbbell are twice that of the
first dumbbell. The initial conditions are chosen such that the
total linear momentum in the inertial frame is zero and the total
energy is positive.
\begin{alignat*}{4}
\bar X_{_0} & = \begin{bmatrix}1&0&0.3\end{bmatrix},&\quad
\bar V_{_0} & = \begin{bmatrix}0&1&0\end{bmatrix},&\quad
\bar \Omega_{1_0} & =
\begin{bmatrix}0&0&9\end{bmatrix},&\quad
R_{_0} & =
I_{3\times 3},\\
\bar x_{2_0} & = \begin{bmatrix}-0.33&0&-0.1\end{bmatrix},&
\bar v_{2_0} & = \begin{bmatrix}0&-0.33&0\end{bmatrix},&
\bar \Omega_{2_0} & =\begin{bmatrix}0&0&0\end{bmatrix},&
R_{2_0} & = I_{3\times 3}.
\end{alignat*}

\begin{figure}[htbp]
    \centering
    \includegraphics[width=0.9\textwidth]{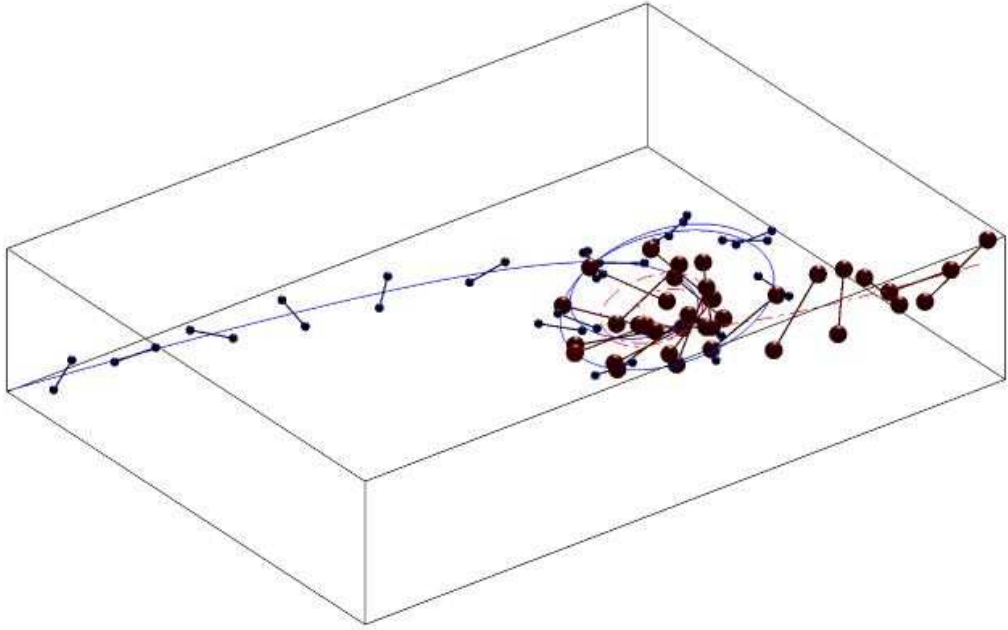}
    \caption{Trajectory in the inertial frame}
    \label{fig:flyby3d}
\end{figure}

\begin{figure}[htbp]
    \centerline{\subfigure[Interchange of energy]{
    \includegraphics[width=0.5\textwidth]{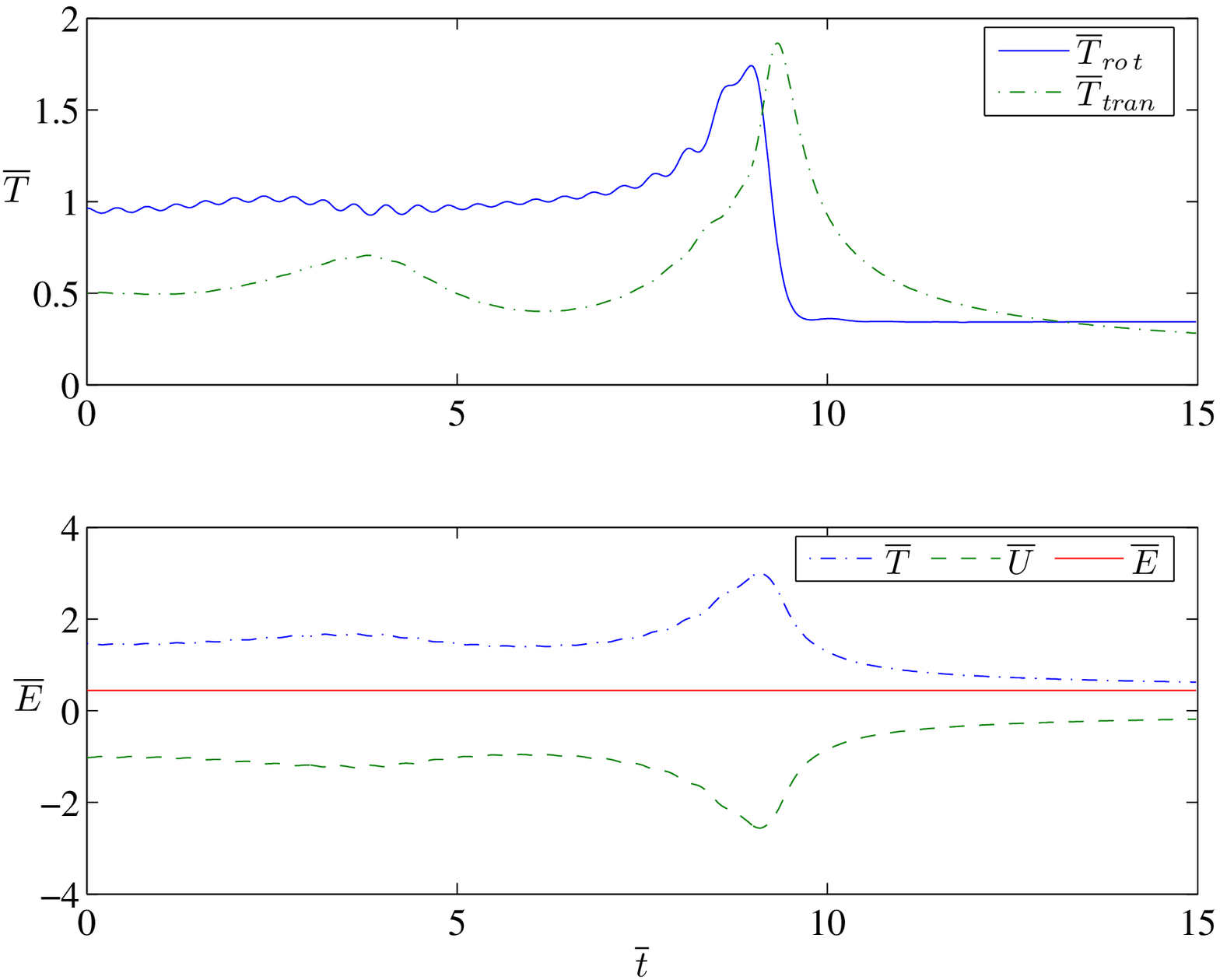}\label{fig:flybyE}}
    \hfill
    \subfigure[Conserved quantities]{
    \includegraphics[width=0.5\textwidth]{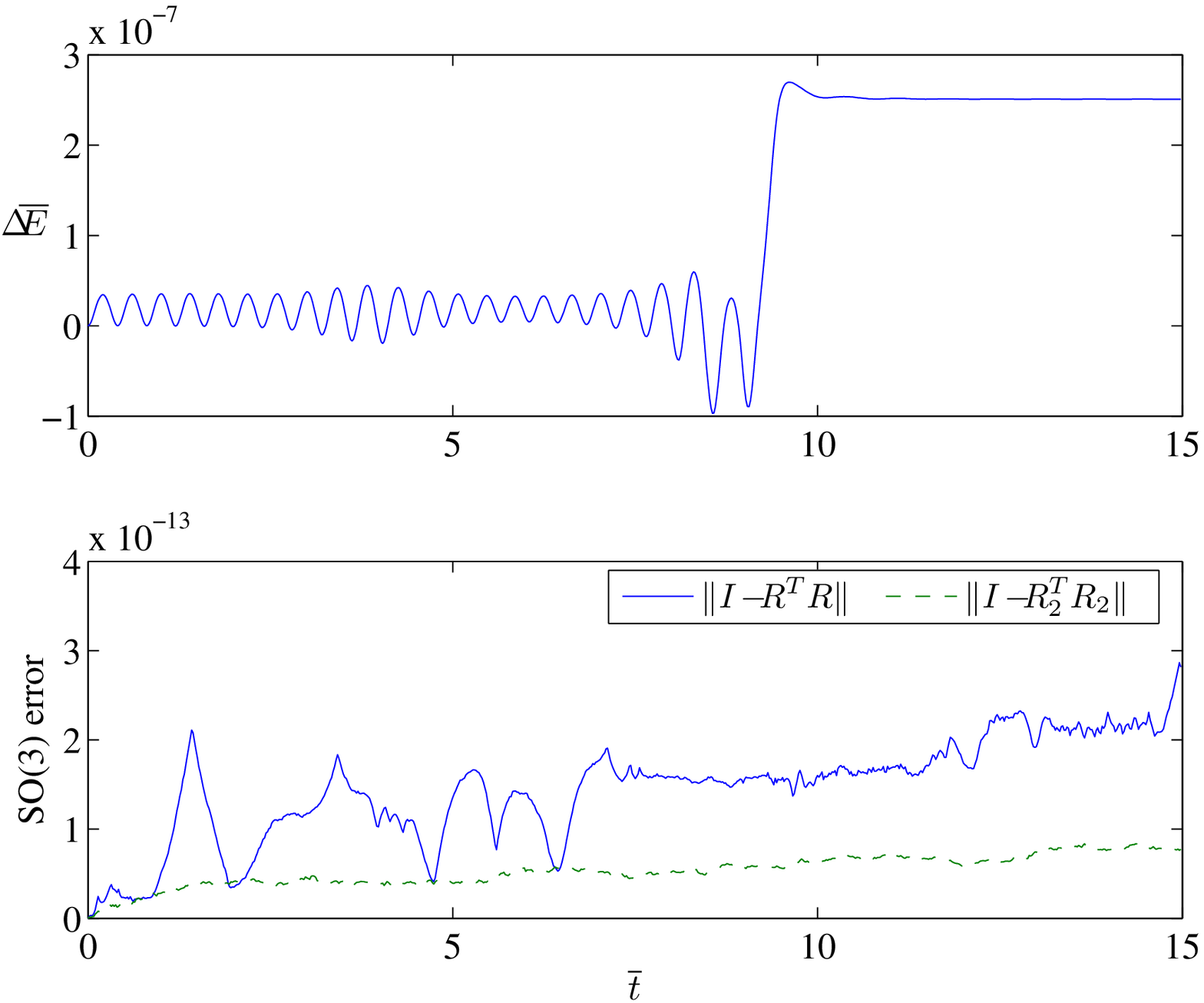}\label{fig:flybyR}}}
    \caption{Lie group variational integrator}
    \label{fig:flybyER}
\end{figure}

Simulation results obtained using the Lie group variational
integrator are given in \reffig{flyby3d} and \reffig{flybyER}.
\reffig{flyby3d} shows the trajectory of the two dumbbells in the
inertial frame. \reffig{flybyE} shows the evolution of the
normalized energy, where the upper figure gives the history of the
translational kinetic energy and the rotational kinetic energy, and
the lower figure shows the interchange between the total kinetic
energy and the gravitational potential energy. \reffig{flybyR} shows
the evolution of the theoretically conserved quantities, where the
upper figure is the history of the total energy, and the lower
figure is the error in the rotation matrix.

Initially, the first dumbbell rotates around the vertical $e_3$
axis, and the second dumbbell does not rotate. Since the angular
velocity of the first dumbbell is relatively large, the rotational
kinetic energy initially exceeds the translational kinetic energy.
As the two dumbbells orbit around each other, the second dumbbell
starts to rotate, the rotational kinetic energy increases, and the
translational kinetic energy decreases slightly for about 6
normalized units of time. At 9 units of time, the distance between
the two dumbbells reaches its minimal separation, and the potential
energy is transformed into kinetic energy, especially translational
kinetic energy. After that, two dumbbells continue to move apart,
and the translational energy and the rotational energy equalize. (A
simple animation of this motion can be found at
\url{http://www.umich.edu/~tylee}.) This shows some of the
interesting dynamics that the full body problem can exhibit. The
non-trivial interchange between rotational kinetic energy,
translational kinetic energy, and potential energy may yield
complicated motions that cannot be observed in the classical two
body problem.

The Lie group variational integrator preserves the total energy and
the geometry of the configuration space. The maximum deviation of
the total energy is $2.6966\times 10^{-7}$, and the maximum value of
the rotation matrix error $\norm{I-R^TR}$ is $2.8657\times
10^{-13}$.

\begin{figure}[tbhp]
    \centerline{\subfigure[Interchange of energy]{
    \includegraphics[width=0.5\textwidth]{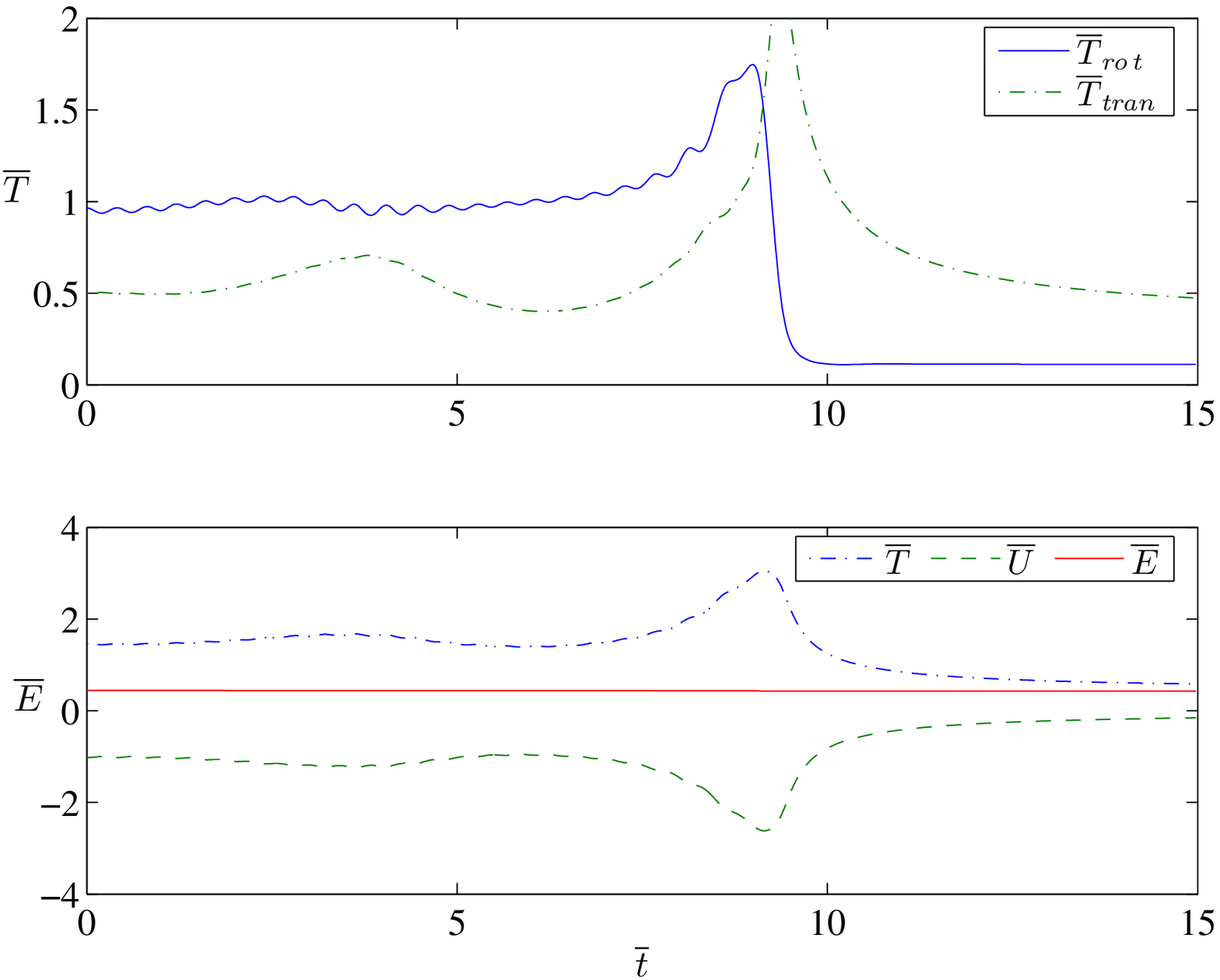}\label{fig:flybyErk}}
    \hfill
    \subfigure[Conserved quantities]{
    \includegraphics[width=0.5\textwidth]{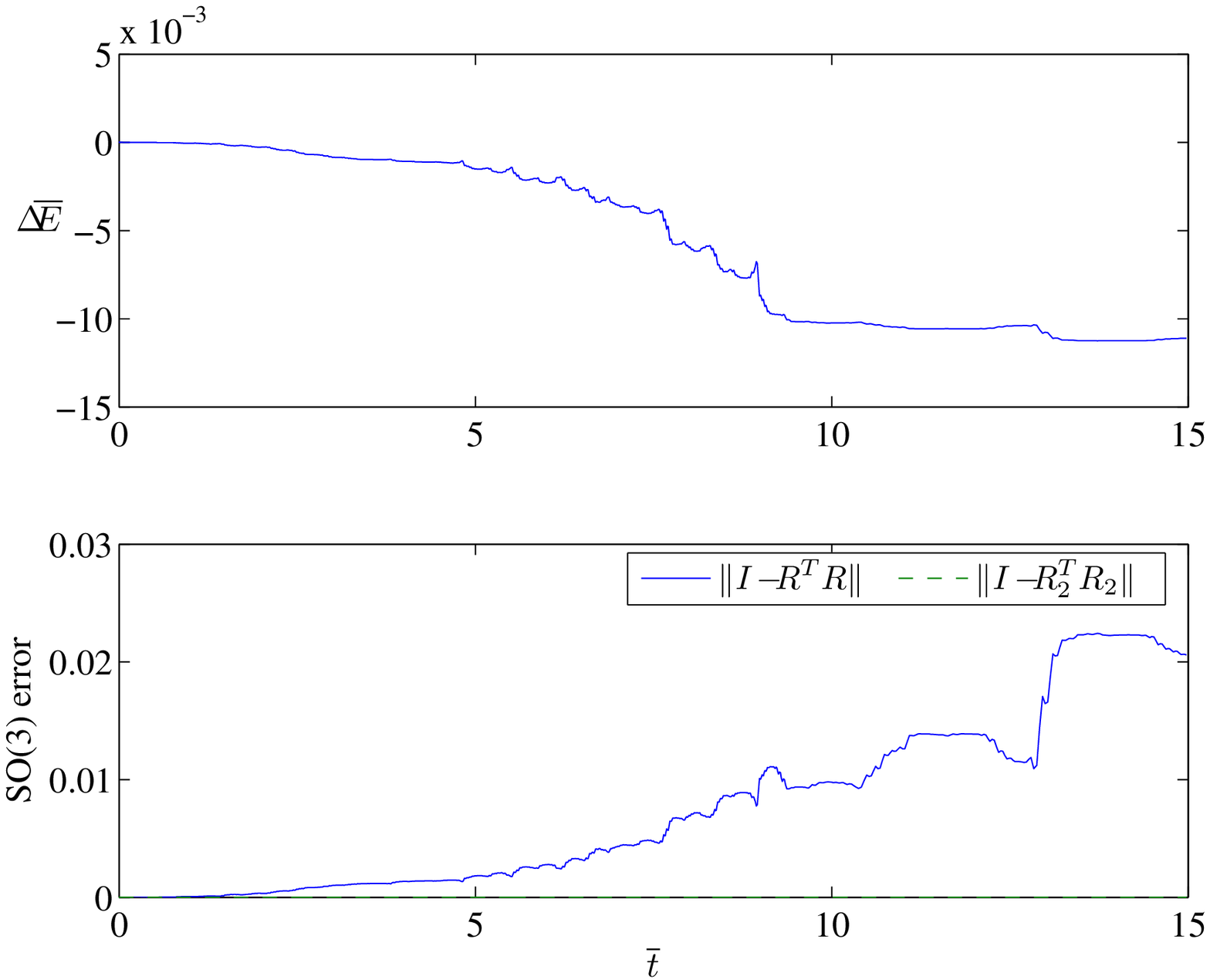}\label{fig:flbyRrk}}}
    \caption{Runge-Kutta method}
    \label{fig:flybyERrk}
\end{figure}

As a comparison, \reffig{flybyERrk} shows simulation results
obtained by numerically integrating the continuous equations of
motion \refeqn{Gamdot}-\refeqn{R2doth} using a standard Runge-Kutta
method. The rotational and the translational kinetic energy
responses are similar to those given in \reffig{flybyER} prior to
the close encounter. However, it fails to simulate the rapid
interchange of the energy near the minimal separation of the two
dumbbells. The deviation of the total energy is relatively large,
with a maximum deviation of $1.1246\times 10^{-2}$. Also, the energy
transfer is quite different from that given in \reffig{flybyE}. The
Runge-Kutta method does not preserve the geometry of the
configuration space, as the discrete trajectory rapidly drifts off
the rotation group to give a maximum rotation matrix error of
$2.2435\times 10^{-2}$. As the gravity and momentum between the two
dumbbells depend on the relative attitude, the errors in the
rotation matrix limits the applicability of standard techniques to
long time simulations.

\section{Conclusions}

\begin{figure}[htbp]
\begin{center}
\resizebox{0.90\linewidth}{!}{
\footnotesize{ 
\xymatrix@!R@!C@C=-0.5cm@R=2cm{ &
\txt{Hamilton\\\refeqn{gamidot}-\refeqn{Ridoth}}
\ar@{-->}[rr]\ar@{.>}'[d][dd] & &
\txt{Discrete Hamilton\\\refeqn{updatexikH}-\refeqn{updateRikH}} \ar@{.>}[dd]^{Reduction}\\
    \txt{Euler-Lagrange\\\refeqn{vidot}-\refeqn{Ridot}} \ar@{-->}[rr]\ar@{->}[ur]\ar@{.>}[dd] & &
\txt{Discrete Euler-Lagrange\\\refeqn{updatexik}-\refeqn{updateRik}}\ar@{->}[ur]\ar@{.>}[dd]\\
    & \txt{Reduced H\\\refeqn{Gamdot}-\refeqn{R2doth}} \ar@{-->}'[r][rr] & & \txt{Reduced DH\\\refeqn{updateX}-\refeqn{updateR2}}\\
    \txt{Reduced EL\\\refeqn{Vdot}-\refeqn{R2dot}} \ar@{-->}[rr]_{Discretization}\ar@{->}[ur]
    & & \txt{Reduced
    DEL\\\refeqn{updateXl}-\refeqn{updatex2l}}\ar@{->}[ur]_{Legendre\;
    trans.}}}
}
\caption{Commutative cube of the equations of
motion}\label{fig:cubic}
\end{center}
\end{figure}
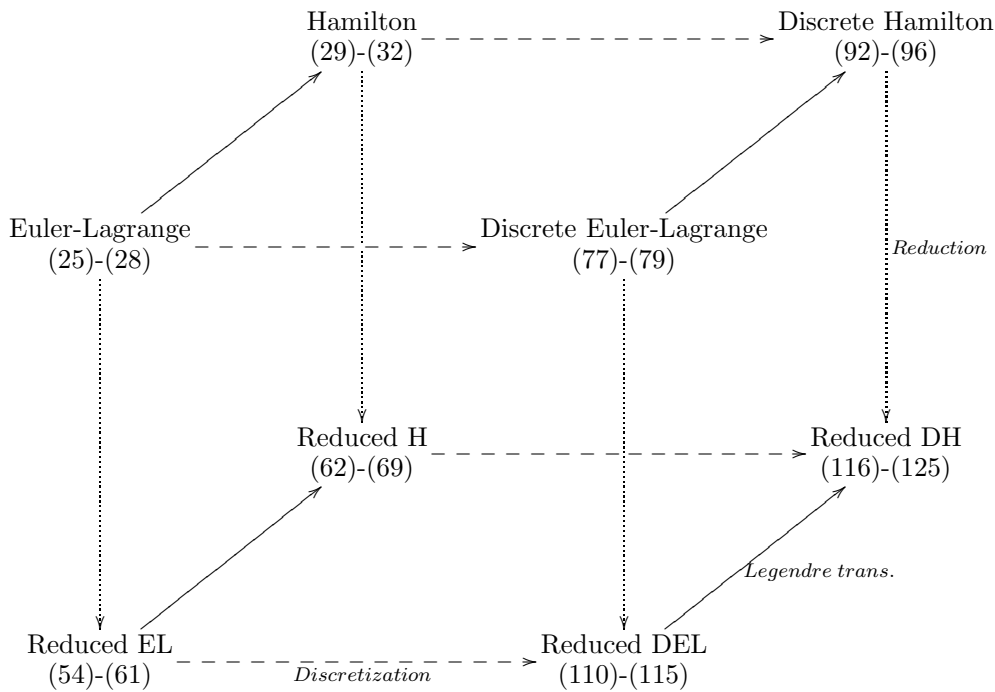

Eight different forms of the equations of motion for the full body
problem are derived. The continuous equations of motion and
variational integrators are derived both in inertial coordinates and
in relative coordinates, and each set of equations of motion is
expressed in both Lagrangian and Hamiltonian form. The relationships
between these equations of motion are summarized in \reffig{cubic}.
This commutative cube was originally given in~\cite{jo:Leok}. In the
figure, dashed arrows represent discretization from the continuous
systems on the left face of the cube to the discrete systems on the
right face. Vertical arrows represent reduction from the full
(inertial) equations on the top face to the reduced (relative)
equations on the bottom face. Front and back faces represent
Lagrangian and Hamiltonian forms, respectively. The corresponding
equation numbers are also indicated in parentheses.

It is shown that the equations of motion for the full body problem
can be derived systematically, using proper Lie group variations,
from Hamilton's principle. The proposed variational integrators
preserve the momenta and symplectic form of the continuous dynamics,
exhibit good energy properties, and they also conserve the geometry
of the configuration space since they are based on Lie group
computations. The main contribution of this paper is the combination
of variational integrators and Lie group computations, developed for
the full body problem. Hence, the resulting numerical integrators
conserve the first integrals as well as the geometry of the
configuration space of the full body dynamics.

\appendix
\section{Appendix}
\subsection{Variations of reduced variables}\label{apprv}
The variations of the reduced variables given in \refeqn{delX}
through \refeqn{delV2} are derived in this section. The variations
of the reduced variables can be obtained from the definitions of the
reduced variables, and the variations of the original variables.

The variation of $X=R_2^T(x_1-x_2)$ is given by
\begin{align*}
\delta X & = \delta R_2^T (x_1-x_2) + R_2 (\delta x_1-\delta x_2).
\end{align*}
Substituting \refeqn{delRi} into the above equation, we obtain
\begin{align*}
\delta X & = - \eta_2 R_2^T (x_1-x_2) + R_2 (\delta x_1-\delta x_2),\\
& = -\eta_2 X + \chi,
\end{align*}
where the reduced variation $\chi:[t_0,t_f]\mapsto\Rset^3$ is
defined to be $\chi=R_2(\delta x_1-\delta x_2)$.

From the definition of $\Omega=R\Omega_1$ and \refeqn{SR},
$S(\delta\Omega)$ is given by
\begin{align*}
S(\delta\Omega)&=\frac{d}{d\epsilon}\bigg|_{\epsilon=0}
S(R^\epsilon
\Omega_1^\epsilon)=\frac{d}{d\epsilon}\bigg|_{\epsilon=0}
R^\epsilon S(\Omega_1^\epsilon)R^{\epsilon T},\\
& = \delta R S(\Omega_1) R^T + R S(\delta \Omega_1) R^T
+ R S(\Omega_1) \delta R^T.
\end{align*}
Substituting \refeqn{delR} and \refeqn{delOmegai} into the above
equation, we obtain
\begin{align*}
S(\delta\Omega)&=\braces{\eta-\eta_2}RS(\Omega_1)R^T
+ R \braces{\dot\eta_1 + S(\Omega_1)\eta_1-\eta_1 S(\Omega_1)}
R^T\\
& \quad + R S(\Omega_1) R^T \braces{-\eta+\eta_2},\\
&=\braces{\eta-\eta_2}S(R\Omega_1)
+ R\dot\eta_1R^T + S(R\Omega_1)R\eta_1R^T-R\eta_1R^T
S(R\Omega_1)\\
& \quad + S(R\Omega_1) \braces{-\eta+\eta_2}.
\end{align*}
Since $\eta=R\eta_1R^T$ and $\Omega=R\Omega_1$, the above equation
reduces to
\begin{align}
S(\delta\Omega) = -\eta_2 S(\Omega) + R\dot{\eta}_1R^T
+S(\Omega)\eta_2.\label{eqn:delOmega0}
\end{align}
From the definition of $R=R_2^TR_1$, $\dot{R}$ is given by
\begin{align}
\dot{R} &= \dot{R}_2^T R_1 + R_2^T \dot{R}_1,\nonumber\\
& = -S(\Omega_2) R + S(\Omega)R.\label{eqn:Rdot0}
\end{align}
Then, $\dot{\eta}$ can be written as
\begin{align}
\dot\eta &= R\dot{\eta}_1R^T + \dot{R}\eta_1R^T +
R\eta_1\dot{R}^T,\nonumber\\
& = R\dot{\eta}_1R^T + \braces{S(\Omega)-S(\Omega_2)}\eta - \eta
\braces{S(\Omega)-S(\Omega_2)}.\label{eqn:etadot}
\end{align}
Substituting \refeqn{etadot} into \refeqn{delOmega0}, we obtain
$S(\delta \Omega)$ in terms of $\eta,\eta_2$ as
\begin{align*}
S(\delta\Omega) & = \dot\eta -S(\Omega)\eta+\eta
S(\Omega)+S(\Omega)\eta_2-\eta_2 S(\Omega)+S(\Omega_2)\eta-\eta
S(\Omega_2),
\end{align*}
which is equivalent to \refeqn{delOmega}.

The variation of $V=R_2^T (\dot{x}_1-\dot{x}_2)$ is given by
\begin{align}
\delta V & = \delta R_2^T (\dot{x}_1-\dot{x}_2)
+ R_2^T(\delta\dot{x}_1-\delta\dot{x}_2),\nonumber\\
& = -\eta_2 V +
R_2^T(\delta\dot{x}_1-\delta\dot{x}_2).\label{eqn:delV0}
\end{align}
From the definition of $\chi=R_2^T(\delta x_1-\delta x_2)$,
$\dot{\chi}$ is given by
\begin{align}
\dot{\chi} & = \dot{R}_2^T(\delta x_1-\delta x_2)
+R_2^T(\delta x_1-\delta x_2),\nonumber\\
& = -S(\Omega_2)\chi+R_2^T(\delta x_1-\delta x_2).\label{eqn:chidot}
\end{align}
Substituting \refeqn{chidot} into \refeqn{delV0}, we obtain
\begin{align*}
\delta V & =-\eta_2 V +\dot\chi + S(\Omega_2)\chi,
\end{align*}
which is equivalent to \refeqn{delV}. The variation $\delta V_2$ can
be derived in the same way, and $S(\delta\Omega_2)$ is given in
\refeqn{delOmegai}.

\subsection{Variations of discrete reduced variables}\label{appdrv}
The variation of the reduced variables $\delta F_{_{k}}$ given in
\refeqn{delFk} is derived in this section. From \refeqn{delFik} and
\refeqn{Fk}, the variation $\delta F_{1_{k}}$ is written as
\begin{align*}
\delta F_{1_{k}} & = -\eta_{1_{k}} F_{1_{k}} + F_{1_{k}} \eta_{1_{k+1}},\\
& = -R_{_{k}}^T \eta_{_{k}}F_{_{k}}R_{_{k}}+R_{_{k}}^T F_{_{k}}
R_{_{k}} R_{_{k+1}}^T \eta_{_{k+1}} R_{_{k+1}},
\end{align*}
where $\eta_{_{k}}\in\so$ is defined as
$\eta_{_{k}}=R_{_{k}}\eta_{1_{k}}R_{_{k}}^T$. Since $F_{_{k}}
R_{_{k}} R_{_{k+1}}^T=F_{_{k}} R_{_{k}}
(R_{_{k}}^TF_{_{k}}^TF_{2_{k}})=F_{2_{k}}$, we have
\begin{align*}
\delta F_{1_{k}} & = R_{_{k}}^T
\parenth{-\eta_{_{k}}F_{_{k}}+F_{2_{k}}
\eta_{_{k+1}}F_{2_{k}}^TF_{_{k}}} R_{_{k}}.
\end{align*}
Then, the variation $\delta F_{_{k}}$ is given by
\begin{align*}
\delta F_{_{k}} & = \delta R_{_{k}} F_{1_{k}} R_{_{k}}^T
+R_{_{k}} \delta F_{1_{k}} R_{_{k}}^T
+R_{_{k}} F_{1_{k}} \delta  R_{_{k}}^T,\nonumber\\
& = -\eta_{2_{k}}F_{_{k}}
+F_{2_{k}}\eta_{_{k+1}}F_{2_{k}}^TF_{_{k}}
+F_{_{k}}\parenth{-\eta_{_{k}} +\eta_{2_{k}} },
\end{align*}
which is equivalent to \refeqn{delFk}.

\bibliography{lgvifbp}
\bibliographystyle{plain}

\end{document}